\documentclass[11pt]{article}
\usepackage{amsmath}
\numberwithin{equation}{subsection}

\usepackage{amsthm}
\usepackage{bbm}
\usepackage{amsmath}
\usepackage{enumerate}
\usepackage{amssymb}
\usepackage{latexsym}
\usepackage{amsfonts}
\usepackage{color}
\usepackage{mathrsfs}
\usepackage{epsfig}
\newtheorem{theorem}{\bf Theorem}[section]
\newtheorem{proposition}[theorem]{\bf Proposition}
\newtheorem{definition}[theorem]{\bf Definition}

\newtheorem{remark}[theorem]{\bf Remark}
\newtheorem{lemma}[theorem]{\bf Lemma}

\newsavebox{\savepar}

\begin{document}
\title{\sc Nodal solutions for  Neumann systems  \\  with gradient dependence} 
\author{\sc 
Kamel Saoudi$^a$,  
Eadah  Alzahrani$^a$,
Du\v{s}an D. Repov\v{s}$^{b,c,d}$\footnote{Corresponding author}\\[1mm]
\small{$^a$ Basic and Applied Scientifc Research Center, Imam Abdulrahman Bin Faisal University,}\\ 
\small{P.O. Box 1982, 31441, Dammam, Saudi Arabia, {\it kmsaoudi@iau.edu.sa, ealzahrani@iau.edu.sa}}\\
\small{$^b$ Faculty of Education, University of Ljubljana, Ljubljana, 1000, Slovenia, {\it dusan.repovs@pef.uni-lj.si}}\\ 
\small{$^c$ Faculty of Mathematics and Physics, University of Ljubljana, Ljubljana, 1000, Slovenia, {\it dusan.repovs@fmf.uni-lj.si}}\\ 		
\small{$^d$ Institute of Mathematics and Physics, Ljubljana, 1000, Slovenia, {\it dusan.repovs@guest.arnes.si}}
} 
		\date{}
		\maketitle		
\begin{abstract}
We consider the following  convective Neumann systems:
\begin{equation*}
\left( \mathrm{S}\right) \qquad \left\{ 
\begin{array}{ll}
-\Delta _{p_1}u_1+\frac{|\nabla u_1|^{p_1}}{u_1+\delta_1 }=f_1(x,u_1,u_2,\nabla u_1,\nabla u_2) & \text{in}%
\;\Omega , \\ 
-\Delta _{p_2}u_2+\frac{|\nabla u_2|^{p_2}}{u_2+\delta_2 }=f_2(x,u_1,u_2,\nabla u_1,\nabla u_2) & \text{in}%
\;\Omega , \\ 
|\nabla u_1|^{p_1-2}\frac{\partial u_1}{\partial \eta }=0=|\nabla u_2|^{p_2-2}\frac{\partial u_2}{\partial \eta } & \text{on}\;\partial
\Omega ,%
\end{array}%
\right.
\end{equation*}
where $\Omega $ is a bounded domain in $\mathbb{R}^{N}$ ($N\geq 2$) with a
smooth boundary $\partial \Omega $, $\delta_1,\,\delta_2 >0$ are small parameters, $\eta $
is the outward unit  vector normal to $\partial \Omega,$
 $f_1,\,f_2:\Omega \times 
\mathbb{R}^2\times \mathbb{R}^{2N}\rightarrow \mathbb{R}$ are 
Carath\'{e}odory
functions that satisfy certain growth conditions,
and
  $\Delta _{p_i}$ ($1<p_i<N,$ for $i=1,2$) 
  are the $p$-Laplace operators
   $\Delta _{p_i}u_i=\mathrm{div}(|\nabla u_i|^{p_i-2}\nabla u_i)$, $\hbox{for every}
\,u_i\in W^{1,p_i}(\Omega ).$ 
In order to prove the existence of solutions to such  systems, we use a sub-supersolution method. We also obtain  nodal solutions
by constructing appropriate sub-solution and super-solution  pairs.  To the best of our knowledge, such systems have not been studied yet.


{\it Math. Subj. Classif. (2020)}:~ 35J62, 35J92.

\end{abstract}
\author{}

\section{Introduction}\label{S1}

In this paper we consider the following  Neumann systems with gradient dependence:
\begin{equation*}
\left( \mathrm{S}\right) \qquad \left\{ 
\begin{array}{ll}
-\Delta _{p_1}u_1+\frac{|\nabla u_1|^{p_1}}{u_1+\delta_1 }=f_1(x,u_1,u_2,\nabla u_1,\nabla u_2) & \text{in}%
\;\Omega , \\ 
-\Delta _{p_2}u_2+\frac{|\nabla u_2|^{p_2}}{u_2+\delta_2 }=f_2(x,u_1,u_2,\nabla u_1,\nabla u_2) & \text{in}%
\;\Omega , \\ 
|\nabla u_1|^{p_1-2}\frac{\partial u_1}{\partial \eta }=0=|\nabla u_2|^{p_2-2}\frac{\partial u_2}{\partial \eta } & \text{on}\;\partial
\Omega ,%
\end{array}%
\right.
\end{equation*}%
where $\Omega $ is a bounded domain in $\mathbb{R}^{N}$ ($N\geq 2$) with a
smooth boundary $\partial \Omega $, $\delta_1,\,\delta_2 >0$ are small parameters, $\eta $
is the outward unit normal vector to $\partial \Omega.$ Here,  $\Delta _{p_i}$ ($1<p_i<N,$ for $i=1,2$) denotes the $p$-Laplace operator,
namely $\Delta _{p_i}u_i:=\mathrm{div}(|\nabla u_i|^{p_i-2}\nabla u_i)$, $\hbox{for every}
\,u_i\in W^{1,p_i}(\Omega ).$   

In recent years,  a lot of work have been done regarding the existence of  solutions
for nonlinear systems with the Dirichlet condition and the reaction term depending on the gradient using different techniques, mainly fixed point theory, variational methods, truncation methods and sub-supersolution methods. We mention for instance,  Candito et al. \cite{CLMous}, where the authors investigated a quasilinear singular Dirichlet system with gradient dependence. They combined Schauder’s fixed
point theorem with sub-supersolutions approach in order to establish the  existence of smooth positive solutions.  For more details, we refer  the readers to 
some recent articles:
Carl and Motreanu    \cite{CaMo},
 Infante et al.   \cite{InMaPr}, 
Miyagaki and Rodrigues     \cite{MiRo}, 
Kita and  Otani \cite{KiOt},
Motreanu et al. \cite{MMM}, 
Orpel \cite{Or},
Ou \cite{Ou}, 
Wang et al. \cite{WaYaGu}, 
Yang and Yang \cite{YaYa}, 
and the references therein. 
See also the
monograph 
by
Motreanu \cite{Mo}. 

On the other hand, the corresponding Neumann system has been much less studied. In this context, the Neumann quasilinear equation involving a connective term  equation was studied in Moussaoui et al. \cite{MoSa}. In Candito et al. \cite{CaMaMo} nodal solutions were obtained for
a $(p_1, p_2)$-Laplacian Neumann system without gradient terms. Neumann systems involving variable exponent double
phase operators and gradient dependence  were investigated in  Guarnotta et al. \cite{GuLiWi}. 

The main interest of the present work  is the presence of the gradient term which constitutes
a serious obstacle in the investigation  of  system $(\mathrm{S})$. 
We note that  system $(\mathrm{S})$  is not in the variational form. Therefore, the usual critical point theory cannot be applied directly. This difficulty is overcome by using the theory of pseudomonotone operators.  We first introduce an auxiliary system using truncation operators. Then we construct a sub-solution $(\underline{u}_1,\underline{u}_2)$ and a super-solution $(\overline{u}_1,\overline{u}_2)$ such that 
$\underline{u}_1\leq \overline{u}_1$, $\underline{u}_2\leq \overline{u}_2$ (see Theorem \ref{thmsubsuper}). At the end, sub and super-solutions and truncation techniques provide  at least two solutions for system $(\mathrm{S}),$
with precise sign properties. 

We shall assume that the nonlinearities $f_i$ for $i=1,2$ are Carath\'{e}odory
functions 
 $f_1,\,f_2:\Omega \times 
\mathbb{R}^2\times \mathbb{R}^{2N}\rightarrow \mathbb{R},$ that
is,  $f_i (., s_1,s_2, \xi_1,\xi_2 )$ is  measurable  for  every  $(s_1,s_2, \xi_1,\xi_2 )\in \mathbb{R}^2\times\mathbb{R}^{2N}$,  and  $f_i (., s_1,s_2, \xi_1,\xi_2 )$  is  continuous  for a.e. $x\in\Omega,$ and they satisfy 
the following growth conditions:
\begin{itemize}
\item[$({\bf H_1})$] There exist $\alpha_i ,\beta_i ,M_i>0$ for $i=1,2,$ such that 
$\max \{\alpha_i ,\beta_i \}<p_i-1$
and 
\begin{equation*}
\Big|f_i(x,s_1,s_2,\xi_1,\xi_2 )\Big|\leq M_i(1+|s_i|^{\alpha_i })(1+|\xi_i |^{\beta_i })\text{,}\ \text{for} \
i=1,2
\ \text{ and  all} \
(x,s_1,s_2,\xi_1,\xi_2 )\in \Omega \times  \mathbb{R}^{2}\times\mathbb{R}^{2N}.
\end{equation*}%
\end{itemize}
\begin{itemize}
\item[$({\bf H_2})$] With appropriate $m_i>0$ for $i=1,2,$ one has 
\begin{equation*}
\underset{|s_i|\rightarrow 0}{\liminf} \Big\{f_i(x,s_1,s_2,\xi_1,\xi_2 ):(\xi_1,\xi_2) \in 
\mathbb{R}^{2N}\Big\}>m_i\text{,} \
\text{uniformly in} \
x\in \Omega.
\end{equation*}%
\end{itemize}
Our main results are the following theorems.

\begin{theorem}\label{T1} Let $\delta_1,\;\delta_2>0$  be small enough and suppose that conditions $({\bf H_1})$
and
$({\bf H_2})$ 
 are satisfied.
  Then system
$\left( \mathrm{S}\right) $ has  a nodal solution $(u_{0},v_0)\in \mathcal{C}
^{1,\gamma }(\overline{\Omega })\times \mathcal{C}
^{1,\gamma }(\overline{\Omega })$ for some
 $\gamma \in (0,1),$ such that $u_{0}(x)$ and $u'_{0}(x)$
are negative near $\partial\Omega.$  
\end{theorem}

\begin{theorem}\label{Thm2}  Let $\delta_1,\;\delta_2>0$  be small enough and suppose that conditions $({\bf H_1})$
and
$({\bf H_2})$ 
 are satisfied.   Then system 
$\left( \mathrm{S}\right) $ has  a positive
solution $(u_{+},u^{+})\in \mathcal{C}^{1,\gamma }(\overline{\Omega })\times \mathcal{C}^{1,\gamma }(\overline{\Omega })$ for some 
 $\gamma \in (0,1),$  such that $u_{+}(x)$ and $u^{+}(x)$ are negative near $\partial\Omega.$  
\end{theorem}

The paper is organized as follows.
 In Section \ref{S3} we collect definitions and results needed in the paper. 
In Section \ref{subsupapproach}, we study auxiliary systems. 
In Section \ref{s4} we prove Theorem~\ref{T5}.
In Section \ref{s5} we study sub-supersolutions.
In Section \ref{Nodal}  we study nodal solutions.
 In Section \ref{s6} we prove our main results. 
 
\section{Preliminaries}\label{S3}

This part is devoted to summarizing the key necessary basic definitions, notations,
and function spaces which will be used in the paper. 
For all other necessary material we refer the reader to the comprehensive monograph
by 
 Papageorgiou et al. \cite{PRR}.
First, some definitions of function spaces. The Banach space $W^{1,p}(\Omega )$ will be equipped with the
usual norm
\begin{equation*}
\Vert u\Vert _{1,p}:=\Big( \Vert u\Vert _{p}^{p}+\Vert \nabla u\Vert
_{p}^{p}\Big)^{1/p}\text{,\quad} \text{for every } u\in W^{1,p}(\Omega )\text{,}
\end{equation*}%
where, 
\begin{equation*}
\Vert v\Vert _{p}:=\left\{ 
\begin{array}{l}
\left( \int_{\Omega }|v(x)|^{p}dx\right) ^{1/p}\text{ if }p<+\infty , \\ 
\\
\text{ess\ }\underset{x\in \Omega }{\sup }\text{\thinspace }|v(x)| \quad \ \ \text{ otherwise.}%
\end{array}%
\right.
\end{equation*}%
Moreover, we
shall
 denote%
\begin{equation*}
\mathcal{W}=W^{1,p_1}(\Omega )\times W^{1,p_2}(\Omega ),\qquad
W_{b}^{1,p_i}(\Omega ):=W^{1,p_i}(\Omega )\cap L^{\infty }(\Omega ),
\end{equation*}
\begin{equation*}
[u_1,u_2]:=\Big\{ u\in W^{1,p}(\Omega ):\; u_1\leq u\leq u_2\Big\},\qquad
 C^{1,\gamma}_0(\overline{\Omega}):=\Big\{ u\in C^{1,\gamma}(\overline{\Omega }): \; u\backslash_{\partial\Omega}=0\Big\}.
\end{equation*}
Now, we define a weak solution of system $(\mathrm{S})$ as follows.
\begin{definition} We say that $(u_1,\,u_2)\in \mathcal{W}$ is a weak solutions of system $(\mathrm{S})$ if
 $$u_i+\delta_i >0
 \hbox{ a.e. in  }
 \Omega,  \quad
 \frac{|\nabla u_i|^{p_i}}{%
u_i+\delta_i }\in L^{1}(\Omega ),
\text{ for }
i=1,2,$$ 
\begin{eqnarray}
&&\int_{\Omega }|\nabla u_1|^{p_1-2}\nabla u_1\nabla \varphi_1 \text{\thinspace }%
dx+\int_{\Omega }\frac{|\nabla u_1|^{p_1}}{u_1+\delta_1 }\varphi_1 \text{\thinspace }%
dx=\int_{\Omega }f_1(x,u_1,u_2,\nabla u_1,\nabla u_2)\varphi_1 \text{\thinspace }dx,\nonumber\\ 
&&\nonumber\\ 
&&\int_{\Omega }|\nabla u_2|^{p_2-2}\nabla u_2\nabla \varphi_2 \text{\thinspace }%
dx+\int_{\Omega }\frac{|\nabla u_2|^{p_2}}{u_2+\delta_2 }\varphi_2 \text{\thinspace }%
dx=\int_{\Omega }f_2(x,u_1,u_2,\nabla u_1,\nabla u_2)\varphi_2 \text{\thinspace }dx,\nonumber\\  \label{defsol}
\end{eqnarray}%
for every $(\varphi_1,\varphi_2) \in W_{b}^{1,p_1}(\Omega )\times W_{b}^{1,p_2}(\Omega ).$
\end{definition}
\begin{remark}
Note that the condition  that $(\varphi_1,\varphi_2) $ are bounded is necessary since $\frac{|\nabla
u_i|^{p}}{u_i+\delta_i }$ for $i=1,2,$ 
is only in $L^{1}(\Omega ).$
\end{remark}

Next, we state the definition of  a sub-solution and a super-solution of system $(\mathrm{S}).$ 
\begin{definition} We say that the pair $(\underline{u}_1,\,\underline{u}_2)\in \mathcal{W}$ is a sub-solution of system $(\mathrm{S})$
 if
 $$\underline{u}_i+\delta_i >0
 \text{ a.e. in }
 \Omega,
 \quad
 \frac{|\nabla \underline{u}_i|^{p_i}}{%
\underline{u}_i+\delta_i }\in L^{1}(\Omega )
\text{ for }
i=1,2,$$ 
\begin{eqnarray}
&&\int_{\Omega }|\nabla \underline{u}_1|^{p_1-2}\nabla \underline{u}_1\nabla \varphi_1 \text{\thinspace }%
dx+\int_{\Omega }\frac{|\nabla \underline{u}_1|^{p_1}}{\underline{u}_1+\delta_1 }\varphi_1 \text{\thinspace }%
dx-\int_{\Omega }f_1(x,\underline{u}_1,w_2,\nabla \underline{u}_1,\nabla w_2)\varphi_1 \text{\thinspace }dx\nonumber\\  
&&+\int_{\Omega }|\nabla \underline{u}_2|^{p_2-2}\nabla \underline{u}_2\nabla \varphi_2 \text{\thinspace }%
dx+\int_{\Omega }\frac{|\nabla \underline{u}_2|^{p_2}}{\underline{u}_2+\delta_2 }\varphi_2 \text{\thinspace }%
dx-\int_{\Omega }f_2(x,w_1,\underline{u}_2,\nabla w_1,\nabla \underline{u}_2)\varphi_2 \text{\thinspace }dx\leq 0,  \nonumber\\\label{defsubsol}
\end{eqnarray}
and we say that the pair $ (\overline{u}_1,\,\overline{u}_2)\in \mathcal{W}$ is a super-solution of
 system $(\mathrm{S})$ if
 $$\overline{u}_i+\delta_i >0
 \text{ a.e. in }
 \Omega,
 \quad
 \frac{|\nabla \overline{u}_i|^{p_i}}{%
\overline{u}_i+\delta_i }\in L^{1}(\Omega )
\text{  for }
i=1,2,$$ 
\begin{eqnarray}
&&\int_{\Omega }|\nabla \overline{u}_1|^{p_1-2}\nabla \overline{u}_1\nabla \varphi_1 \text{\thinspace }%
dx+\int_{\Omega }\frac{|\nabla \overline{u}_1|^{p_1}}{\overline{u}_1+\delta_1 }\varphi_1 \text{\thinspace }%
dx-\int_{\Omega }f_1(x,\overline{u}_1,w_2,\nabla \overline{u}_1,\nabla w_2)\varphi_1 \text{\thinspace }dx\nonumber\\  
&&+\int_{\Omega }|\nabla \overline{u}_2|^{p_2-2}\nabla \overline{u}_2\nabla \varphi_2 \text{\thinspace }%
dx+\int_{\Omega }\frac{|\nabla \overline{u}_2|^{p_2}}{\overline{u}_2+\delta_2 }\varphi_2 \text{\thinspace }%
dx-\int_{\Omega }f_2(x,w_1,\overline{u}_2,\nabla w_1,\nabla \overline{u}_2)\varphi_2 \text{\thinspace }dx\geq 0,\nonumber\\  \label{defsupersol}
\end{eqnarray}%
for every $(\varphi_1,\varphi_2) \in W_{b}^{1,p_1}(\Omega )\times W_{b}^{1,p_2}(\Omega )$ with $\varphi_1,\,\varphi_2 \geq 0$ in $\Omega$ and  for all $(w_1,w_2)\in\mathcal{W},$ such that $\underline{u}_i\leq w_i\leq \overline{u}_i,$ for $i=1,2$ and  with  all  integrals  in   \eqref{defsubsol}  and   \eqref{defsupersol}  being  ﬁnite.
\end{definition}
We shall be using the following conditions in the paper: 
\begin{itemize}
\item[$({\bf H_3})$] Let $0\leq q_1\leq p_1-1$ and $0\leq r_1\leq p_2-1$. For every $\rho >0,$ there
exists $M_1:=M_1(\rho )>0$ such that 
\begin{equation*}
|f_1(x,s_1,s_2,\xi_1,\xi_2 )|\leq M_1(1+|\xi_1 |^{q_1}+|\xi_2 |^{r_1})\text{ \ in }\Omega \times \lbrack -\rho
,\rho ]^2\times \mathbb{R}^{2N}.
\end{equation*}
\item[$({\bf H_4})$] Let $0\leq q_2\leq p_1-1$ and $0\leq r_2\leq p_2-1$. For every $\rho >0,$ there
exists $M_2:=M_2(\rho )>0$ such that 
\begin{equation*}
|f_2(x,s_1,s_2,\xi_1,\xi_2 )|\leq M_2(1+|\xi_1 |^{q_2}+|\xi_2 |^{r_2})\text{ \ in }\Omega \times \lbrack -\rho
,\rho ]^2\times \mathbb{R}^{2N}.
\end{equation*}
\item[$({\bf H_5})$] There are $\underline{u}_1,\overline{u}_1\in \mathcal{C}%
^{1}(\overline{\Omega })$ sub- and supersolution of system $(\mathrm{S})$ respectively, satisfying
\begin{equation}
\overline{u}_1+\delta_1 \geq \underline{u}_1+\delta_1 >0\text{ \ a.e. in }\Omega ,
\label{c3}
\end{equation}%
\item[$({\bf H_6})$] There are $\underline{u}_2,\overline{u}_2\in \mathcal{C}%
^{1}(\overline{\Omega })$ sub and supersolution of system $(\mathrm{S})$ respectively, satisfying
\begin{equation}
\overline{u}_2+\delta_2 \geq \underline{u}_2+\delta_2 >0\text{ \ a.e. in }\Omega ,
\label{c4}
\end{equation}
\end{itemize}

Via a standard argument, we shall  prove the following.
\begin{proposition}\label{lemma-1}
Suppose  that conditions  $({\bf H_3})$, $({\bf H_4})$, $({\bf H_5}),$ and $({\bf H_6})$ are satisfied. Let
$(\underline{u}_i,\underline{v}_i),\, (\overline{u}_i,\overline{v}_i) \in W^{1,p_1}_b(\Omega )\times W^{1,p_2}_b(\Omega )$ be  pairs of sub-solution and super-solutions for system $\left( \mathrm{S}\right).$  Set $$ \overline{u}=\min \{\overline{u}_{1},%
\overline{u}_{2}\}\qquad \underline{u}=\max \{\underline{u}_{1},%
\underline{u}_{2}\}$$
$$ \overline{v}=\min \{\overline{v}_{1},%
\overline{v}_{2}\}\qquad \underline{v}=\max \{\underline{v}_{1},%
\underline{v}_{2}\}$$
and assume that $ \overline{u}\leq \overline{v}$ and $\underline{u}\leq \underline{v}$ Then $(\underline{u},\underline{v}),(\overline{u},\overline{v})$  is also a pair of sub-solution and  super-solution for system 
$\left( \mathrm{S}\right).$
\end{proposition}

\begin{proof} The proof was inspired by the proof of Motreanu et al. \cite[ Lemma 3]{MMP}. Fix $\epsilon >0$ and  define the truncation function $%
\xi _{\epsilon }(s)=\max \{-\epsilon ,\min \{s,\epsilon \}\},$ for every $s\in 
\mathbb{R}.$ By Marcus et al. \cite{MM}, we know that  
$$\xi _{\epsilon }((\overline{u}_{1}-%
\overline{u}_{2})^{-}),\; \xi _{\epsilon }((\underline{u}_{1}-%
\underline{u}_{2})^{+})\in \mathcal{W},$$%
\begin{equation*}
\nabla \xi _{\epsilon }((\overline{u}_{1}-\overline{u}_{2})^{-})=\xi
_{\epsilon }^{\prime }((\overline{u}_{1}-\overline{u}_{2})^{-})\nabla (%
\overline{u}_{1}-\overline{u}_{2})^{-}
\end{equation*}
and 
\begin{equation*}
\nabla \xi _{\epsilon }((\underline{u}_{1}-\underline{u}_{2})^{+})=\xi
_{\epsilon }^{\prime }((\underline{u}_{1}-\underline{u}_{2})^{+})\nabla (%
\underline{u}_{1}-\underline{u}_{2})^{+}.
\end{equation*}
Now, letting $\varphi \in C_{c}^{1}(\Omega )$ be a test function with $\varphi \geq 0,$ we obtain
\begin{eqnarray}\label{eq20}
\left\langle -\Delta_{{p}_1}\underline{u}_{1}+\frac{|\nabla \underline{u}_{1}|^{p_1}}{\underline{u}_{1}+\delta_1},\xi _{\epsilon }((\underline{u}_{1}-\underline{u}_{2})^{+})\varphi \right\rangle \nonumber\\
\leq \int_{\Omega }f_1(x,\underline{u}_{1},w_2,\nabla \underline{u}_{1},\nabla w_2 )\xi
_{\epsilon }((\underline{u}_{1}-\underline{u}_{2})^{+})\varphi \text{ }%
\mathrm{d}x,%
\end{eqnarray}
\begin{eqnarray}\label{eq21}
\left\langle -\Delta_{{p}_1}\overline{u}_{1}+\frac{|\nabla \overline{u}_{1}|^{p_1}%
}{\overline{u}_{1}+\delta_1 },\xi _{\epsilon }((\overline{u}_{1}-\overline{u%
}_{2})^{-})\varphi \right\rangle \nonumber\\
\geq \int_{\Omega }f_1(x,\overline{u}_{1},w_2,\nabla \overline{u}_{1},\nabla w_2)\xi
_{\epsilon }((\overline{u}_{1}-\overline{u}_{2})^{-})\varphi \text{ }%
\mathrm{d}x,%
\end{eqnarray}
for every $w_2\in W^{1,p_2}(\Omega)$ with $\underline{u}_{2}\leq w_2\leq \overline{u}_{2}$  and%

\begin{eqnarray}\label{eq22}
\left\langle -\Delta_{{p}_1}\underline{u}_{2}+\frac{|\nabla \underline{u}_{2}|^{p_1}%
}{\underline{u}_{2}+\delta_1 },(\epsilon -\xi _{\epsilon }((\underline{u}%
_{1}-\underline{u}_{2})^{+}))\varphi \right\rangle \nonumber\\
\leq \int_{\Omega }f_1(x,w_1,\underline{u}_{2},\nabla w_1, \nabla\underline{u}_{2})\left(
\epsilon -\xi _{\epsilon }((\underline{u}_{1}-\underline{u}%
_{2})^{+})\right) \varphi \text{ }\mathrm{d}x,
\end{eqnarray}

\begin{eqnarray}\label{eq23}
\left\langle -\Delta_{{p}_1}\overline{u}_{2}+\frac{|\nabla \overline{u}_{2}|^{p_1}%
}{\overline{u}_{2}+\delta_1 },(\epsilon -\xi _{\epsilon }((\overline{u}%
_{1}-\overline{u}_{2})^{-}))\varphi \right\rangle \nonumber\\ 
\geq \int_{\Omega }f_1(x,w_1,\overline{u}_{2},\nabla w_1, \nabla\overline{u}_{2})\left(
\epsilon -\xi _{\epsilon }((\overline{u}_{1}-\overline{u}%
_{2})^{-})\right) \varphi \text{ }\mathrm{d}x,
\end{eqnarray}
for every $w_1\in W^{1,p_1}(\Omega)$ with $\underline{u}_{1}\leq w_1\leq \overline{u}_{1}.$   Therefore, by  the monotonicity of the 
$-p$-Laplacian operator, one has
\begin{eqnarray}\label{eq24}
&&\\
&&\left\langle -\Delta_{{p}_1}\underline{u}_{1}+\frac{|\nabla \underline{u}_{1}|^{p_1}}{\underline{u}_{1}+\delta_1},\xi _{\epsilon }((\underline{u}_{1}-\underline{u}_{2})^{+})\varphi \right\rangle+ \left\langle -\Delta_{{p}_1}\underline{u}_{2}+\frac{|\nabla \underline{u}_{2}|^{p_1}%
}{\underline{u}_{2}+\delta_1 },(\epsilon -\xi _{\epsilon }((\underline{u}%
_{1}-\underline{u}_{2})^{+}))\varphi \right\rangle \nonumber\\
&&\geq \int_{\Omega }|\nabla \underline{u}_{1}|^{p_1-2}(\nabla \underline{u}%
_{1},\nabla \varphi )_{
\mathbb{R}
^{N}}\xi _{\epsilon }((\underline{u}_{1}-\underline{u}_{2})^{+})\text{ }%
\mathrm{d}x+\int_{\Omega }\frac{|\nabla \underline{u}_{1}|^{p_1}}{\underline{u}%
_{1}+\delta_1 }\xi _{\epsilon }((\underline{u}_{1}-\underline{u}%
_{2})^{+})\varphi \text{ }\mathrm{d}x \nonumber\\  
&&+\int_{\Omega }|\nabla \underline{u}_{2}|^{p_1-2}(\nabla \underline{u}%
_{2},\nabla \varphi )_{
\mathbb{R}
^{N}}\left( \epsilon -\xi _{\epsilon }((\underline{u}_{1}-\underline{u}%
_{2})^{+})\right) \text{ }\mathrm{d}x+\int_{\Omega }\frac{|\nabla \underline{u%
}_{2}|^{p_1}}{\underline{u}_{2}+\delta_1 }(\epsilon -\xi _{\epsilon }((%
\underline{u}_{1}-\underline{u}_{2})^{+}))\varphi \text{ }\mathrm{d}x.\nonumber
\end{eqnarray}
and 
\begin{eqnarray}\label{eq25}
&&\\
&&\left\langle -\Delta_{{p}_1}\overline{u}_{1}+\frac{|\nabla \overline{u}_{1}|^{p_1}}{\overline{u}_{1}+\delta_1},\xi _{\epsilon }((\overline{u}_{1}-\overline{u}_{2})^{-})\varphi \right\rangle+ \left\langle -\Delta_{{p}_1}\overline{u}_{2}+\frac{|\nabla \overline{u}_{2}|^{p_1}%
}{\overline{u}_{2}+\delta_1 },(\epsilon -\xi _{\epsilon }((\overline{u}%
_{1}-\overline{u}_{2})^{-}))\varphi \right\rangle \nonumber\\
&&\leq \int_{\Omega }|\nabla \overline{u}_{1}|^{p_1-2}(\nabla \overline{u}%
_{1},\nabla \varphi )_{
\mathbb{R}
^{N}}\xi _{\epsilon }((\overline{u}_{1}-\overline{u}_{2})^{-})\text{ }%
\mathrm{d}x+\int_{\Omega }\frac{|\nabla \overline{u}_{1}|^{p_1}}{\overline{u}%
_{1}+\delta_1 }\xi _{\epsilon }((\overline{u}_{1}-\overline{u}%
_{2})^{-})\varphi \text{ }\mathrm{d}x \nonumber\\  
&&+\int_{\Omega }|\nabla \overline{u}_{2}|^{p_1-2}(\nabla \overline{u}%
_{2},\nabla \varphi )_{
\mathbb{R}
^{N}}\left( \epsilon -\xi _{\epsilon }((\overline{u}_{1}-\overline{u}%
_{2})^{+})\right) \text{ }\mathrm{d}x+\int_{\Omega }\frac{|\nabla \overline{u%
}_{2}|^{p_1}}{\overline{u}_{2}+\delta_1 }(\epsilon -\xi _{\epsilon }((%
\overline{u}_{1}-\overline{u}_{2})^{-}))\varphi \text{ }\mathrm{d}x.\nonumber
\end{eqnarray}
Invoking equations \eqref{eq20}, \eqref{eq22}, and \eqref{eq24}, we obtain
\begin{eqnarray*}
&&\int_{\Omega }|\nabla \underline{u}_{1}|^{p_1-2}(\nabla \underline{u}_{1},\nabla
\varphi )_{
\mathbb{R}
^{N}}\frac{1}{\epsilon }\xi _{\epsilon }((\underline{u}_{1}-\underline{u}%
_{2})^{+})\text{ }\mathrm{d}x+\int_{\Omega }\frac{|\nabla \underline{u}%
_{1}|^{p_1}}{\underline{u}_{1}+\delta_1 }\frac{1}{\epsilon }\xi _{\epsilon
}((\underline{u}_{1}-\underline{u}_{2})^{-}\text{ }\mathrm{d}x \\  
&&+\int_{\Omega }|\nabla \underline{u}_{2}|^{p_1-2}(\nabla \underline{u}%
_{2},\nabla \varphi )_{
\mathbb{R}
^{N}}\left( 1-\frac{1}{\epsilon }\xi _{\epsilon }((\underline{u}_{1}-%
\underline{u}_{2})^{+})\right) \text{ }\mathrm{d}x\\
&&+\int_{\Omega }\frac{%
|\nabla \underline{u}_{2}|^{p_1}}{\underline{u}_{2}+\delta_1 }(1-\frac{1}{%
\epsilon }\xi _{\epsilon }((\underline{u}_{1}-\underline{u}_{2})^{+})%
\text{ }\mathrm{d}x \\ 
&&\leq \int_{\Omega }f_1(x,\underline{u}_{1},w_2,\nabla \underline{u}_{1},\nabla w_2)\frac{1}{%
\epsilon }\xi _{\epsilon }((\underline{u}_{1}-\underline{u}%
_{2})^{+})\varphi \text{ }\mathrm{d}x\\
&&+\int_{\Omega }f_1(x,\underline{u}%
_{1},w_2,\nabla \underline{u}_{1},\nabla w_{2})\left( 1-\frac{1}{\epsilon }\xi
_{\epsilon }((\underline{u}_{1}-\underline{u}_{2})^{-})\right) \varphi 
\text{ }\mathrm{d}x.%
\end{eqnarray*}
In a similar manner,  invoking equations \eqref{eq21}, \eqref{eq23}, and  \eqref{eq25}, we get
\begin{eqnarray*}
&&\int_{\Omega }|\nabla \overline{u}_{1}|^{p_1-2}(\nabla \overline{u}_{1},\nabla
\varphi )_{%
\mathbb{R}
^{N}}\frac{1}{\epsilon }\xi _{\epsilon }((\overline{u}_{1}-\overline{u}%
_{2})^{-})\text{ }\mathrm{d}x+\int_{\Omega }\frac{|\nabla \overline{u}%
_{1}|^{p_1}}{\overline{u}_{1}+\delta_1 }\frac{1}{\epsilon }\xi _{\epsilon
}((\overline{u}_{1}-\overline{u}_{2})^{-}\text{ }\mathrm{d}x \\ 
&&+\int_{\Omega }|\nabla \overline{u}_{2}|^{p-2}(\nabla \overline{u}%
_{2},\nabla \varphi )_{
\mathbb{R}
^{N}}\left( 1-\frac{1}{\epsilon }\xi _{\epsilon }((\overline{u}_{1}-%
\overline{u}_{2})^{-})\right) \text{ }\mathrm{d}x\\
&&+\int_{\Omega }\frac{%
|\nabla \overline{u}_{2}|^{p_1}}{\overline{u}_{2}+\delta_1 }(1-\frac{1}{%
\varepsilon }\xi _{\epsilon }((\overline{u}_{1}-\overline{u}_{2})^{-})%
\text{ }\mathrm{d}x \\  
&&\geq \int_{\Omega }f_2(x,w_{1},\overline{u}_{2},\nabla w_{1},\nabla \overline{u}_{2})\frac{1}{%
\epsilon }\xi _{\epsilon }((\overline{u}_{1}-\overline{u}%
_{2})^{-})\varphi \text{ }\mathrm{d}x\\
&&+\int_{\Omega }f_2(x,w_{1},\overline{u}%
_{2},\nabla w_{1},\nabla \overline{u}_{2})\left( 1-\frac{1}{\epsilon }\xi
_{\epsilon }((\overline{u}_{1}-\overline{u}_{2})^{-})\right) \varphi 
\text{ }\mathrm{d}x.%
\end{eqnarray*}
Letting $\epsilon \rightarrow 0,$  and observing that
$$
\begin{cases}
\frac{1}{\epsilon }\xi _{\epsilon }((\overline{u}_{1}-\overline{u}_{2})^{-}\rightarrow \mathrm{1}_{\{\overline{u}_{1}<\overline{u}_{2}\}}(x)%
\text{, \ a.e. in }\Omega \text{ as }\epsilon \rightarrow 0,\\
\\
\frac{1}{\epsilon }\xi _{\epsilon }((\underline{u}_{1}-\underline{u}%
_{2})^{+}\rightarrow \mathrm{1}_{\{\underline{u}_{1}<\underline{u}_{2}\}}(x)%
\text{, \ a.e. in }\Omega \text{ as }\epsilon \rightarrow 0,
\end{cases}
$$
we see that 
\begin{eqnarray*}
&&\int_{\Omega }|\nabla \underline{u}|^{p_1-2}\nabla \underline{u}\nabla \varphi 
\text{ }\mathrm{d}x+\int_{\Omega }\frac{|\nabla \underline{u}|^{p_1}}{\underline{%
u}+\delta_1 }\varphi \text{ }\mathrm{d}x\leq \int_{\Omega }f_1(x,\underline{u}_1,w_2
,\nabla\underline{u}_1,w_2)\varphi \text{ }\mathrm{d}x,\\
&&\textrm{and}\\
&&\int_{\Omega }|\nabla \overline{u}|^{p_1-2}\nabla \overline{u}\nabla \varphi 
\text{ }\mathrm{d}x+\int_{\Omega }\frac{|\nabla \overline{u}|^{p_1}}{\overline{%
u}+\delta_1 }\varphi \text{ }\mathrm{d}x\geq \int_{\Omega }f_1(x,\overline{u}_1,w_2
,\nabla\overline{u}_1,\nabla w_2)\varphi \text{ }\mathrm{d}x,
\end{eqnarray*}%
for every $\varphi \in C_{c}^{1}(\Omega ),$ $\varphi \geq 0$ a.e. in $\Omega.$
 By a similar argument as above, we obtain
\begin{eqnarray*}
&&\int_{\Omega }|\nabla \underline{v}|^{p_2-2}\nabla \underline{v}\nabla \varphi 
\text{ }\mathrm{d}x+\int_{\Omega }\frac{|\nabla \underline{v}|^{p_2}}{\underline{%
v}+\delta_2 }\varphi \text{ }\mathrm{d}x\leq \int_{\Omega }f_2(x,w_1,\underline{v}_2,\nabla w_1,\nabla\underline{v}_2)\varphi \text{ }\mathrm{d}x,\\
&&\textrm{and}\\
&&\int_{\Omega }|\nabla \overline{v}|^{p_2-2}\nabla \overline{v}\nabla \varphi 
\text{ }\mathrm{d}x+\int_{\Omega }\frac{|\nabla \overline{v}|^{p_2}}{\overline{%
v}+\delta_2 }\varphi \text{ }\mathrm{d}x\geq \int_{\Omega }f_2(x,w_1,\overline{v}_2,\nabla w_1,\nabla\overline{v}_2)\varphi \text{ }\mathrm{d}x.
\end{eqnarray*}%
Finally, in view of the density of   $C_{c}^{1}(\Omega )$  in both $W^{1,p_1}(\Omega )$ and $W^{1,p_2}(\Omega ),$ we can deduce that $(\underline{u},\underline{v}),(\overline{u},\overline{v})$  is also a pair of sub-solution and  super-solution of system 
$\left( \mathrm{S}\right).$ 
\end{proof}

\section{Auxiliary systems}\label{subsupapproach}

Let, the pair  $(\underline{u}_1,\underline{u}_2),\;(\overline{u}_1,\overline{u}_2)\in \mathcal{C}^{1}(\overline{\Omega })\times \mathcal{C}^{1}(\overline{\Omega })$ are 
sub-solutions and super-solutions respectively of system $(\mathrm{S})$ as required in conditions $({\bf H_5})$ and $({\bf H_6})$. Now, for a given $(u_1,u_2)\in \mathcal{W},$ 
we shall define 
the truncation operators $\mathcal{T}_i%
:W^{1,p_i}(\Omega )\rightarrow W^{1,p_i}(\Omega )$  by%
\begin{eqnarray}
\mathcal{T}_1(u_1):=\begin{cases}
\underline{u}_1 & \text{when }u_1\leq \underline{u}_1, \\ 
u_1 & \text{if }\underline{u}_1\leq u_1\leq \overline{u}_1, \\ 
\overline{u}_1 & \text{otherwise,}%
\end{cases}
\qquad\qquad \mathcal{T}_2(u_2):=
\begin{cases}
\underline{u}_2 & \text{when }u_2\leq \underline{u}_2, \\ 
u_2 & \text{if }\underline{u}_2\leq u_2\leq \overline{u}_2, \\ 
\overline{u}_2 & \text{otherwise.}%
\end{cases}%
\nonumber\\ \label{eq1}
\end{eqnarray}
Then 
by Carl et al. \cite[Lemma 2.89]{CLM}, 
 $\mathcal{T}_1,$ $\mathcal{T}_2$ are continuous, monotone and bounded. In view of conditions $({\bf H_3})$ and  $({\bf H_4}),$ if $\rho >0$ then 
\begin{equation}
-\rho \leq \underline{u}_1\leq \overline{u}_1\leq \rho,\qquad\qquad -\rho \leq \underline{u}_2\leq \overline{u}_2\leq \rho.  \label{eq2}
\end{equation}%
We introduce the Nemitskii operators $\mathcal{N}_{{f}_1}$ and $\mathcal{N}_{{f}_2}$ generated by the
Carath\'{e}odory functions $f_1$ and $f_2$ respectively, which are well defined for $i=1,2$ since the range of   
$\mathcal{T}_1,\; \mathcal{T}_2$ lies within the region $[\underline{u}_i, \overline{u}_i].$ So, due to the Rellich-Kondrachov Compactness Embedding Theorem, the maps
\begin{equation}\label{eq3}
\mathcal{N}_{{f}_1}\circ(\mathcal{T}_1,\mathcal{T}_2):[\underline{u}_1,\overline{u}_1]\subset
\mathcal{W} \longrightarrow L^{p_1^{\prime }}(\Omega )\hookrightarrow W^{-1,p_1}(\Omega ),
\end{equation}
\begin{equation}\label{eq4}
\mathcal{N}_{{f}_2}\circ(\mathcal{T}_1,\mathcal{T}_2):[\underline{u}_2,\overline{u}_2]\subset
\mathcal{W} \longrightarrow L^{p_2^{\prime }}(\Omega )\hookrightarrow W^{-1,p_2}(\Omega )
\end{equation}
  are bounded and completely continuous. Furthermore, set
  $$\mathcal{F}(u)=\Big(\mathcal{N}_{{f}_1}(\mathcal{T}_1u_1, \mathcal{T}_2u_2,\nabla(\mathcal{T}_1u_1),\nabla(\mathcal{T}_2u_2),\mathcal{N}_{{f}_2}(\mathcal{T}_1u_1, \mathcal{T}_2u_2,\nabla(\mathcal{T}_1u_1),\nabla(\mathcal{T}_2u_2)\Big).$$
    Next,  define the cut-off functions $b_i:\Omega \times \mathbb{R}%
\longrightarrow \mathbb{R}$ for $i=1,2$  by
\begin{equation}\label{eq5}
b_1(x,s):=-(\underline{u}_1(x)-s)_{+}^{p_1-1}+(s-\overline{u}_1(x))_{+}^{p_1-1},\quad \text{ for every } 
(x,s)\in \Omega \times \mathbb{R},
\end{equation}
\begin{equation}\label{eq6}
b_2(x,s):=-(\underline{u}_2(x)-s)_{+}^{p_2-1}+(s-\overline{u}_2(x))_{+}^{p_2-1},\quad  \text{ for every } 
(x,s)\in \Omega \times \mathbb{R}.
\end{equation}
It is easy
 to see  that 
 $b_i$
 are Carath\'{e}odory functions
  for $i=1,2,$ 
   fulfilling the following growth
condition 
\begin{equation}\label{growth}
|b_1(x,s)|\leq \varphi_1(x)+c_1|s|^{p_1-1},\ \text{for a.e. }x\in \Omega \text{ and every }%
s\in 
\mathbb{R},
\end{equation}
\begin{equation}\label{growth1}
|b_2(x,s)|\leq \varphi_2(x)+c_2|s|^{p_2-1},\ \text{for a.e. }x\in \Omega  \ \text{and  every }%
s\in 
\mathbb{R},
\end{equation}
with  $\varphi_1,\;\varphi_2$ $\in L^{\infty }(\Omega )$ and $c_1,\,c_2> 0.$ 
Moreover, one has the following estimates
\begin{equation}\label{eq7}
\int_{\Omega }b_1(\cdot ,u_1)u_1\,\mathrm{d}x\geq C_{1}\Vert u_1\Vert _{p_1}^{p}-C_{2},%
\text{ for every }u_1\in W^{1,p_1}(\Omega ), 
\end{equation}%
\begin{equation}\label{eq8}
\int_{\Omega }b_2(\cdot ,u_2)u_2\,\mathrm{d}x\geq C'_{1}\Vert u_1\Vert _{p_1}^{p}-C'_{2},%
\text{ for every }u_1\in W^{1,p_2}(\Omega ), 
\end{equation}%
where  $C_{1},\;C_{2},\;C'_{1},\;C'_{2}$ are some positive constants (for more details see, e.g., Carl et al. \cite[pp. 95--96]{CLM}). Let $\mu> 0$  and  set 
$$\mu \mathcal{B}(u)=\Big(\mu\mathcal{B}_1(u_1),\mu\mathcal{B}_2(u_2)\Big).$$
Now, we introduce the following   auxiliary problem:
\begin{equation*}
\left( \mathrm{S}_\mu\right) \qquad \left\{ 
\begin{array}{ll}
-\Delta _{p_1}u_1+\frac{|\nabla(\mathcal{T} u_1)|^{p_1}}{\mathcal{T}u_1+\delta_1 }=f_1(x,\mathcal{T}u_1,\mathcal{T}u_2,\nabla (\mathcal{T}u_1),\nabla (\mathcal{T}u_2))-\mu b_1(x,u) & \text{in}%
\;\Omega , \\ 
\\
-\Delta _{p_2}u_2+\frac{|\nabla (\mathcal{T}u_2)|^{p_2}}{\mathcal{T}u_2+\delta_2 }=f_2(x,\mathcal{T}u_1,\mathcal{T}u_2,\nabla (\mathcal{T}u_1),\nabla (\mathcal{T}u_2))-\mu b_2(x,u) & \text{in}%
\;\Omega , \\ 
\\
|\nabla u_1|^{p_1-2}\frac{\partial u_1}{\partial \eta }=0=|\nabla u_2|^{p_2-2}\frac{\partial u_2}{\partial \eta } & \text{on}\;\partial
\Omega ,%
\end{array}%
\right.
\end{equation*}
where $(u_1,u_2)\in \mathcal{W}.$ Our main result in this section  concerning system $(\mathrm{S}_{\mu })$ is as
follows.
\begin{theorem}\label{T5}
Suppose that conditions $({\bf H_3})$, $({\bf H_4})$, $({\bf H_5}),$ and $({\bf H_6})$ are satisfied.
 Then system $(\mathrm{S}_{\mu })$ has a pair of weak solutions $(u_1, u_2)\in \mathcal{W}.$
\end{theorem}
The following estimates will be key for the proof of Theorem \ref{T5} in the next section.
\begin{lemma}\label{L3} Suppose that conditions $({\bf H_3})$ and  $({\bf H_4})$  are satisfied. Then there exist  constants $k_{0}, k'_{0}>0$
 depending only on $p_1,\,p_2,$ and $\Omega$ such that
\begin{eqnarray*}
\int_{\Omega }\Big|f_1\Big(x ,\mathcal{T}u_1,\mathcal{T}u_2,\nabla (\mathcal{T}u_1),\nabla (\mathcal{T}u_2)\Big)\Big||u_1|\mathrm{d}%
x&\leq& \frac{1}{2}\Big(\Vert \nabla u_1\Vert _{p_1}^{p_1}+\Vert \nabla u_2\Vert _{p_2}^{p_2}\Big)\\
&+&k_{0}\Big(1+\left\Vert
u_1\right\Vert _{p_1}+\left\Vert u_1\right\Vert _{p_1}^{p_1}+\left\Vert u_1\right\Vert _{p_2}^{p_2}\Big).
\end{eqnarray*}
and 
\begin{eqnarray*}
\int_{\Omega }\Big|f_2\Big(x ,\mathcal{T}u_1,\mathcal{T}u_2,\nabla (\mathcal{T}u_1),\nabla (\mathcal{T}u_2)\Big)\Big||u_2|\mathrm{d}%
x&\leq& \frac{1}{2}\Big(\Vert \nabla u_1\Vert _{p_1}^{p_1}+\Vert \nabla u_2\Vert _{p_2}^{p_2}\Big)\\
&+&k'_{0}\Big(1+\left\Vert
u_2\right\Vert _{p_2}+\left\Vert u_2\right\Vert _{p_2}^{p_2}+\left\Vert u_2\right\Vert _{p_1}^{p_1}\Big).
\end{eqnarray*}
for every $(u_1,u_2)\in \mathcal{W}.$
\end{lemma}

\begin{proof} We 
shall
prove only the first inequality. The second inequality can be verified similarly.
First,  by condition $({\bf H_3})$, we have
\begin{eqnarray}\label{eq9}
&&\displaystyle\int_{\Omega }\Big|f_1\Big(x ,\mathcal{T}u_1,\mathcal{T}u_2,\nabla (\mathcal{T}u_1),\nabla (\mathcal{T}u_2)\Big)\Big||u_1|\mathrm{d}x\nonumber\\
&&\leq M_1\int_\Omega\Big(1+|\nabla (\mathcal{T}u_1)|^{q_1}+|\nabla (\mathcal{T}u_2)|^{r_1}\Big)|u_1|\mathrm{d}x. 
\end{eqnarray}%
Now, using Young's inequality, we get for any fixed $\varepsilon \in ]0,\frac{1}{2M}[$ and  every $u_1\in W^{1,p_1}(\Omega )$, that 
\begin{equation}\label{eq10}
|\nabla (\mathcal{T}u_1)|^{q_1}|u_1|\leq \varepsilon |\nabla (\mathcal{T}u_1)|^{%
\frac{q_1p_1}{p_1-1}}+c_{\varepsilon }|u_1|^{p_1}\leq \varepsilon \Big(1+|\nabla (\mathcal{T}u_1)|^{p_1}\Big)+c_{\varepsilon }|u_1|^{p_1}.
\end{equation}
Similarly, for every $u_2\in W^{1,p_2}(\Omega ),$ we have
\begin{equation}\label{eq11}
|\nabla (\mathcal{T}u_2)|^{r_1}|u_1|\leq \varepsilon |\nabla (\mathcal{T}u_2)|^{%
\frac{r_1p_2}{p_2-1}}+c'_{\varepsilon }|u_1|^{p_2}\leq \varepsilon \Big(1+|\nabla (\mathcal{T}u_2)|^{p_2}\Big)+c'_{\varepsilon }|u_1|^{p_2}.
\end{equation}%
On the other hand, using equation \eqref{eq1}, we can see that 
\begin{eqnarray}\label{eq12}
\int_{\Omega }|\nabla (\mathcal{T}u_1)|^{p_1}\text{ }\mathrm{d}x&=&\int_{\{%
\underline{u}_1\leq u_1\leq \overline{u}_1\}}|\nabla u_1|^{p_1}\text{ }\mathrm{d}%
x+\int_{\{u_1\geq \overline{u}_1\}}|\nabla \overline{u}_1|^{p_1}\text{ }\mathrm{d}x+\int_{\{u_1\leq 
\underline{u}_1\}}|\nabla \underline{u}_1|^{p_1}\text{ }\mathrm{d}x
\nonumber\\ 
&\leq&\int_{\Omega }|\nabla u_1|^{p_1}\text{ }\mathrm{d}x+\int_{\Omega }|\nabla \underline{u}_1|^{p_1}\text{ }\mathrm{d}%
x+\int_{\Omega }|\nabla 
\overline{u}_1|^{p_1}\text{ }\mathrm{d}x \nonumber\\ 
&\leq& \Vert \nabla
u_1\Vert _{p_1}^{p_1}+ \Big(\left\Vert \nabla \underline{u}_1\right\Vert _{\infty
}^{p_1}+\left\Vert \nabla \overline{u}_1\right\Vert _{\infty }^{p_1}\Big)|\Omega |.\nonumber\\
\end{eqnarray}
Using the same techniques, we can get
\begin{eqnarray}\label{eq13}
\int_{\Omega }|\nabla (\mathcal{T}u_2)|^{p_2}\text{ }\mathrm{d}x&\leq& \Vert \nabla
u_2\Vert _{p_2}^{p_2}+ \Big(\left\Vert \nabla \underline{u}_2\right\Vert _{\infty
}^{p_2}+\left\Vert \nabla \overline{u}_2\right\Vert _{\infty }^{p_2}\Big)|\Omega |.\nonumber\\
\end{eqnarray}
Consequently, using equations \eqref{eq10}-\eqref{eq13}, ones gets
\begin{eqnarray}
&&\displaystyle\int_{\Omega }\Big|f_1\Big(x ,\mathcal{T}u_1,\mathcal{T}u_2,\nabla (\mathcal{T}u_1),\nabla (\mathcal{T}u_2)\Big)\Big||u_1|\mathrm{d}x\leq
 M_1\Big(|\Omega|^{\frac{p_1-1}{p_1}}||u_1||_{p_1}\nonumber\\
&& +\epsilon|\Omega|\Big(1+||\nabla \underline{u}_1||_{\infty
}^{p_1}+||\nabla \overline{u}_1||_{\infty }^{p_1}\Big)+\epsilon
||\nabla u_1\Vert _{p_1}^{p_1}+c_{\epsilon }|| u_1||_{p_1}^{p_1}\nonumber \\&&+\epsilon|\Omega|\Big(1+||\nabla \underline{u}_2||_{\infty
}^{p_2}+||\nabla \overline{u}_2||_{\infty }^{p_2}\Big)+\epsilon
||\nabla u_2||_{p_2}^{p_2}+c_{\epsilon }|| u_1||_{p_2}^{p_2}\Big)\nonumber \\ 
 &&\leq \frac{1}{2}\Big(\Vert \nabla u_1\Vert _{p_1}^{p_1}+\Vert \nabla u_2\Vert _{p_2}^{p_2}\Big)+k_{0}\Big(1+\left\Vert
u_1\right\Vert _{p_1}+\left\Vert u_1\right\Vert _{p_1}^{p_1}+\left\Vert u_1\right\Vert _{p_2}^{p_2}\Big).\nonumber\\
\end{eqnarray}%
for a suitable $k_{0}>0.$   The proof of  Lemma \ref{L3} is thus completed.
\end{proof}
The following useful  estimates can be verified in a similar way as in Moussaoui et al. \cite[ Lemma 2.2]{MoSa}.
\begin{lemma}\label{L4}Suppose that conditions $({\bf H_3})$, $({\bf H_4})$, $({\bf H_5}),$ and $({\bf H_6})$ are satisfied. Then for
every $u=(u_1,u_2)\in \mathcal{W}$, there exist  constants $k_1,\,k_2$ independent of $u_1,\; u_2$ respectively, such that
\begin{eqnarray}\label{eq14}
&&\frac{|\nabla (\mathcal{T}u_1)|^{p_1}}{%
\mathcal{T}u_1+\delta_1 }|u_1|\in L^{1}(\Omega ) \quad\textrm{and}\quad\displaystyle\int_{\Omega }\frac{|\nabla (\mathcal{T}u_1)|^{p_1}}{\mathcal{T}u_1+\delta_1 }|u_1|%
\text{ }\mathrm{d}x\leq k_1\Big(1+|| u_1||_{p_1}\Big),\nonumber\\
&&
\end{eqnarray}
\begin{eqnarray}\label{eq15}
&&\frac{|\nabla (\mathcal{T}u_2)|^{p_2}}{%
\mathcal{T}u_2+\delta_2 }|u_2|\in L^{1}(\Omega) \quad\textrm{and}\quad 
\displaystyle\int_{\Omega }\frac{|\nabla (\mathcal{T}u_2)|^{p_2}}{\mathcal{T}u_2+\delta_1 }|u_2|%
\text{ }\mathrm{d}x\leq k_2\Big(1+|| u_2||_{p_2}\Big).\nonumber\\
&&
\end{eqnarray}
\end{lemma}

\section{Proof of Theorem~\ref{T5}}\label{s4}

First, by the growth conditions in  equations \eqref{growth}
 and \eqref{growth1} we know that the Nemytskii operators $\mathcal{B}_i:W^{1,p_i}(\Omega
)\longrightarrow W^{-1,p^{\prime }_i}(\Omega )$ given by $\mathcal{B}_i%
u_i(x)=b(\cdot ,u_i)$ are  well defined, continuous and bounded for $i=1,2.$ Also, the operator
$\mathcal{B}(u)=(\mathcal{B}_1(u_1),\mathcal{B}_2(u_2))$ is well defined.  Moreover, using
the compact embedding $W^{1,p_i}(\Omega )\hookrightarrow L^{p_i}(\Omega )$, we know that the operator 
$\mathcal{B}$ is completely continuous. Next, using conditions $(\mathrm{H.}3)$ and $(\mathrm{H.}4)$, we can
 introduce the function $\pi _{p_i,\delta_i }:(-\delta_i ,+\infty )\times \mathbb{R}^{N}\longrightarrow \mathbb{R}$  for $i=1,2$ defined by \begin{equation*}
\pi _{p_i,\delta_i }(s_i,\xi_i )=\frac{|\xi_i |^{p_i}}{s_i+\delta_i }
\end{equation*}%
having the growth
\begin{equation*}
|\pi _{p_i,\delta_i }(s_i,\xi_i )|\leq \delta _{0}|\xi_i |^{p_i}
\end{equation*}
 for every $s_i>-\delta_i$ 
 and all $\xi_i \in 
\mathbb{R}^{N},$ where $\delta _{0}>0$ be a constant such that 
$
\overline{u}_i+\delta_i \geq \underline{u}_i+\delta_i >\delta _{0}$  a.e. in $%
\Omega $ for $i=1,2.$ 

By virtue of Motreanu et al. \cite[ 
Theorem 2.76]{MMPA} and Gasinski et al. \cite[Theorem 3.4.4]{GP}),  we know that the corresponding Nemytskii operator 
$$\Pi _{p_i,\delta_i }:[\underline{u}_i,\overline{u}_i]\subset W^{1,p_i}(\Omega
)\longrightarrow L^{1}(\Omega )\subset W^{-1,p_i^{\prime }}(\Omega )$$ 
is bounded and continuous for $i=1,2.$ By virtue of the compact embedding of $W^{1,p}(\Omega )$
into $L^{p}(\Omega )$, we know that $\Pi_{p,\delta }(u) =(\Pi_{p_1,\delta_1 }(u_1),\Pi_{p_2,\delta_2 }(u_2) )$ 
is completely continuous. Finally,  $\mathcal{A}(u)=(A_1(u_1), A_2(u_2)),$ where $\mathcal{A}: \mathcal{W} \rightarrow \mathcal{W}^*$ defined in equation \eqref{defsol}, is well defined, bounded, continuous, strictly monotone, and of type $(S_+).$ Therefore, for every $u$ and $\varphi\in \mathcal{W},$ we have the following representations
\begin{eqnarray*}
\displaystyle \langle\mathcal{A}(u), \varphi\rangle_{\mathcal{W}}&=&\sum_{i=1}^2\int_{\Omega }|\nabla u_i|^{p_i-2}\nabla u_i\nabla \varphi_i \text{\thinspace }dx,\\
\langle \Pi_{p,\delta }(u),\varphi\rangle_{\mathcal{W}}&=&\sum_{i=1}^2\int_{\Omega } \Pi _{p_i,\delta_i }(u_1,u_2,\nabla u_1,\nabla u_2)\varphi_i dx,\\
\langle \mathcal{B}(u),\varphi\rangle_{\mathcal{W}}&=&\sum_{i=1}^2\int_{\Omega } \mathcal{B}_i(u_1,u_2,\nabla u_1,\nabla u_2)\varphi_i dx,\\
\langle\mathcal{F}(u),\varphi\rangle_{\mathcal{W}}&=&\sum_{i=1}^2\int_{\Omega }\mathcal{N}_{{f}_i}(\mathcal{T}_1u_1, \mathcal{T}_2u_2,\nabla(\mathcal{T}_1u_1),\nabla(\mathcal{T}_2u_2)\varphi_i dx.
\end{eqnarray*}
Now, for every $u$ and $\varphi\in \mathcal{W},$ system $(\mathrm{P}_{\mu })$ can be given in the form
\begin{equation}\label{operator1}
\langle\mathcal{A}(u)+\mu \mathcal{B}u+\Pi _{p,\delta }(u),\varphi\rangle_{\mathcal{W}}=\langle\mathcal{F}(u),\varphi\rangle_{\mathcal{W}}.  
\end{equation}%
Set $$\mathcal{\chi}_\mu:=\mathcal{A}(u)+\mu \mathcal{B}u+\Pi _{p,\delta }(u)-\mathcal{F}(u).$$
First, by conditions  $(\mathrm{H.}1)$ and $(\mathrm{H.}2)$,   $\mathcal{\chi}_{\mu}$ is well defined, continuous and bounded. The next step in the proof is to show that the operator $\mathcal{\chi}_{\mu }$
is pseudo-monotone. To this end,  using the $(\mathrm{S})_{+}$-property of $\mathcal{A}$ and in  view of the compactness of the operators $\Pi_{p,\delta },\;\mathcal{B},\;\mathcal{F}$, we can use Gambera et al. \cite[ Lemma 2.2]{GaGu} to deduce that the operator $\mathcal{\chi}_\mu$ also has the $(\mathrm{S})_{+}$-property. Furthermore, we can apply Zeidler
\cite[Proposition 26.2]{Ze} to see that the operator $\mathcal{\chi}_\mu$ is pseudo-monotone. 

Let us show that the operator $\mathcal{\chi}_{\mu }:\mathcal{W}\rightarrow \mathcal{W}^*$ is coercive. To this end, using equation (\ref{operator1}), we get
\begin{eqnarray}\label{eq16}
\left\langle \mathcal{\chi}_{\mu }(u),u\right\rangle &=&\sum_{i=1}^2\int_{\Omega }|\nabla
u_i|^{p_i}\text{ }\mathrm{d}x+\mu\sum_{i=1}^2 \int_{\Omega }b_i(x,u_i)u_i\text{ }\mathrm{d}
x+\sum_{i=1}^2\int_{\Omega }\frac{|\nabla (\mathcal{T}u_i)|^{p_i}}{\mathcal{T}u_i+\delta_i }u_i
\text{ }\mathrm{d}x\nonumber\\
&-&\sum_{i=1}^2\int_{\Omega }f_i(x ,\mathcal{T}u_1,\mathcal{T}u_2,\nabla (\mathcal{T}u_1),\nabla (\mathcal{T}%
u_2))u_i\,\mathrm{d}x \nonumber\\ 
&\geq&\sum_{i=1}^2\int_{\Omega }|\nabla
u_i|^{p_i}\text{ }\mathrm{d}x+\mu\sum_{i=1}^2 \int_{\Omega }b_i(x,u_i)u_i\text{ }\mathrm{d}
x-\sum_{i=1}^2\int_{\Omega }\frac{|\nabla (\mathcal{T}u_i)|^{p_i}}{\mathcal{T}u_i+\delta_i }u_i
\text{ }\mathrm{d}x\nonumber\\
&-&\sum_{i=1}^2\int_{\Omega }f_i(x ,\mathcal{T}u_1,\mathcal{T}u_2,\nabla (\mathcal{T}u_1),\nabla (\mathcal{T}%
u_2))u_i\,\mathrm{d}x.\nonumber\\
\end{eqnarray}%
Now, using equation \eqref{eq16} and combining equations \eqref{eq7} and \eqref{eq8}, with Lemmas \ref{L3} and \ref{L4}, we obtain
\begin{eqnarray}\label{eq17}
\left\langle \mathcal{\chi}_{\mu }(u),u\right\rangle &\geq&  \sum_{i=1}^2 ||\nabla u_i||
_{p_i}^{p_i}+\mu\Big( C_{1}|| u_1|| _{p_1}^{p_1}-C_{2}\Big)+\mu \Big( C'_{1}|| u_2|| _{p_2}^{p_2}-C'_{2}\Big) \nonumber\\ 
&-&\sum_{i=1}^2 k_i\Big(1+|| u_i|| _{p_i}\Big)- \Big(\Vert \nabla u_1\Vert _{p_1}^{p_1}+\Vert \nabla u_2\Vert _{p_2}^{p_2}\Big)\nonumber\\
&-&k_{0}\Big(1+||u_1||_{p_1}+|| u_1||_{p_1}^{p_1}+||u_1||_{p_2}^{p_2}\Big)\nonumber\\
&-&k'_{0}\Big(1+\left\Vert
u_2\right\Vert _{p_2}+\left\Vert u_2\right\Vert _{p_2}^{p_2}+\left\Vert u_2\right\Vert _{p_1}^{p_1}\Big)\nonumber\\
&\geq&   \Big(\Vert \nabla u_1\Vert _{p_1}^{p_1}+\Vert \nabla u_2\Vert _{p_2}^{p_2}\Big)+\mu C^*_{1}\Big( || u_1|| _{p_1}^{p_1}+|| u_2|| _{p_2}^{p_2}\Big)-\mu \Big( C_{2}+C'_{2}\Big) \nonumber\\ 
&-&\sum_{i=1}^2 k_i\Big(1+|| u_i|| _{p_i}\Big)-k^*_{0}\Big(1+||u_1||_{p_1}+|| u_1||_{p_1}^{p_1}+||u_1||_{p_2}^{p_2}\Big),\nonumber\\
\end{eqnarray}%
where $k^*_{0}:=\max\{k_{0},\;k'_{0}\}$ and $C^*_{1}:=\min\{C_{1},\;C'_{1}\}.$
Then, 
invoking Moussaoui and Saoudi
  \cite[Lemma 2.2]{MoSa}, we can deduce that
$$\displaystyle|| \nabla (u_i\mathbbm{1}_{\{\underline{u}_i<u_i<\overline{u}_i\}})||_{p_i}\leq 
\tilde{C}_i,
\
\hbox{for some}
\
\tilde{C}_i>0,
\
\hbox{independent of}
\
u_i,
\
\hbox{for}
\
i=1,2.$$
 Furthermore, for  sufficiently large
  $\mu >0$
 such that
  $\mu C^*_{1}-k^*_{0}>0,$ 
  and for every sequence 
  $(u_{n})_{n}$ in $\mathcal{W}$, 
   inequality \eqref{eq17} forces 
\begin{equation*}
\langle \mathcal{\chi}_{\mu}(u_{n}),u_{n}\rangle\rightarrow +\infty,\;\;\textrm{as}\;\; ||u_n||_{\mathcal{W}}\rightarrow +\infty.
\end{equation*}

Therefore, since $\mathcal{\chi}_\mu$ is continuous, bounded, coercive and pseudomonotone, invoking
 the 
Pseudomonotone Operators Theorem  (see e. g., Carl et al. \cite[Theorem 2.99]{CLM}), we get the existence of $u\in
\mathcal{W}$ such that 
\begin{equation}\label{eq18}\left\langle \mathcal{\chi}_{\mu }(u_1,u_2),(\varphi_1,\varphi_2)\right\rangle =0,\;\;\textrm{ for every}\;\; (\varphi_1,\varphi_2) \in
\mathcal{W}.
\end{equation}
 Moreover, using Casas et al. \cite[Theorem 3]{CF}, we have
\begin{equation*}
|\nabla u_1|^{p_1-2}\frac{\partial u_1}{\partial \eta }=0=|\nabla u_2|^{p_2-2}\frac{\partial u_2}{\partial \eta }=0\ \text{on }\partial
\Omega \text{.}
\end{equation*}%
Therefore, we can conclude that  $u=(u_1,u_2)\in \mathcal{W}$ is a weak solution of $(\mathrm{S}_{\mu })$.
This completes the proof of Theorem~\ref{T5}.
\qed
\section{Sub-supersolutions}\label{s5}
The aim of this section is to construct  pairs of sub-solution and super-solution of system $(\mathrm{S}).$
\begin{theorem}\label{thmsubsuper}
\label{T3} Assume that conditions $({\bf H_3})$, $({\bf H_4})$, $({\bf H_5}),$ and $({\bf H_6})$ are satisfied. Then
system $(\mathrm{S})$ has a solution $u=(u_1,u_2)\in \mathcal{C}^{1,\gamma }(
\overline{\Omega })\times \mathcal{C}^{1,\gamma }(%
\overline{\Omega })\cap [\underline{u}_1,\overline{u}_1]\times [\underline{u}_2,\overline{u}_2]$ for some $\gamma \in (0,1).$ 
\end{theorem}

\begin{proof}[Proof of Theorem~\ref{thmsubsuper}] First, using Theorem \ref{T5}, we can fix $\mu>0$ sufficiently large such that system $(\mathrm{S}_{\mu })$ admits a pair of weak
solutions $u=(u_1,u_2)\in \mathcal{W}.$
 It remains to verify that $ u=(u_1,u_2)\in [\underline{u}_1,\overline{u}_1]\times [\underline{u}_2,\overline{u}_2].$ Here, we give just the proof for $u_1\in [\underline{u}_1,\overline{u}_1].$ A  similar reasoning  yields
 the second inequality. To this end, we set $(\varphi_1,\varphi_2)=((u_1-\overline{u}_1)_+,0).$ By Lemma \ref{L4} and  condition $({\bf H_5})$ combined with equation \eqref{eq18}, we obtain
 \begin{eqnarray*}
&&\int_{\Omega }|\nabla u_1|^{p_1-2}\nabla u_1\text{\thinspace }\nabla (u_1-\overline{u}_1)_{+}\text{ }\mathrm{d}x+\int_{\Omega }\frac{|\nabla (\mathcal{T}u_1)|^{p_1}}{%
\mathcal{T}u_1+\delta_1 }(u_1-\overline{u}_1)_{+}\text{ }\mathrm{d}x\\ 
&&=\int_{\Omega }f(x ,\mathcal{T}u_1,\mathcal{T}u_2,\nabla (\mathcal{T}u_1),\nabla (\mathcal{T}u_2))(u_1-\overline{u}_1)_{+}\text{ }\mathrm{d}x-\mu \int_{\Omega }b(x ,u_1)(u_1-\overline{u}_1)_{+}
\text{ }\mathrm{d}x \\ 
&&=\int_{\Omega }f(x ,\overline{u}_1,,\mathcal{T}u_2,\nabla \overline{u}_1,\nabla (\mathcal{T}u_2))(u_1-\overline{u}_1)_{+}%
\text{ }\mathrm{d}x-\mu \int_{\Omega }(u_1-\overline{u}_1)_{+}^{p_1}\text{ }%
\mathrm{d}x \\ 
&&\leq \int_{\Omega }|\nabla \overline{u}_1|^{p_1-2}\nabla \overline{u}_1\,\nabla (u_1-\overline{u}_1)_{+}\,\mathrm{d}x+\int_{\Omega }\frac{|\nabla \overline{u}_1|^{p_1}%
}{\overline{u}_1+\delta_1 }(u_1-\overline{u}_1)_{+}\,\mathrm{d}x-\mu \int_{\Omega
}(u_1-\overline{u}_1)_{+}^{p_1}\text{ }\mathrm{d}x\text{.}%
\end{eqnarray*}
Now, according to equation \eqref{eq1},  
\begin{equation*}
\int_{\Omega }\frac{|\nabla (\mathcal{T}u_1)|^{p_1}}{\mathcal{T}u_1+\delta_1 }(u_1-\overline{u}_1)_{+}\text{ }\mathrm{d}x=\int_{\Omega }\frac{|\nabla \overline{u}_1
|^{p_1}}{\overline{u}_1+\delta_1 }(u_1-\overline{u}_1)_{+}\,\mathrm{d}x,
\end{equation*}%
so it follows that 
\begin{equation}\label{eq19}
\int_{\Omega }\Big( |\nabla u_1|^{p_1-2}\nabla u_1-|\nabla \overline{u}_1
|^{p_1-2}\nabla \overline{u}_1\Big) \nabla (u_1-\overline{u}_1)_{+}\,\mathrm{d}%
x\leq -\mu \int_{\Omega }(u_1-\overline{u}_1)_{+}^{p_1}\text{ }\mathrm{d}x\leq 0.
\end{equation}
Hence,  it follows
from equation \eqref{eq19} combined with  the monotonicity of $ A_{1},$ that $u_1\leq \overline{u}_1.$ In
the
 same way, in order to see that 
$\underline{u}_1\leq u_1,$ we set $(\varphi_1,\varphi_2)=((\underline{u}_1-u_1)_+,0).$
 So, $u=(u_1,u_2)\in[\underline{u}_1,\overline{u}_1]\times [\underline{u}_2,\overline{u}_2].$ Moreover, 
 according to Miyajima et al.  \cite[Remark 8 ]{MMT}, we obtain that 
$u=(u_1,u_2)\in \mathcal{C}^{1,\gamma }(
\overline{\Omega })\times \mathcal{C}^{1,\gamma }(%
\overline{\Omega }),$ for some $\gamma \in (0,1)$ and $\frac{\partial u_1}{%
\partial \eta }=\frac{\partial u_2}{%
\partial \eta }=0$ on $\partial \Omega $. Therefore, we have shown that $u=(u_1,u_2)\in \mathcal{C}^{1,\gamma }(
\overline{\Omega })\times \mathcal{C}^{1,\gamma }(%
\overline{\Omega })$ is a solution of the system $(\mathrm{S})$ within $[\underline{u}_1,\overline{u}_1]\times [\underline{u}_2,\overline{u}_2].$
\end{proof}
\section{Nodal solutions}\label{Nodal}
The objective of this section is to show the existence of nodal solutions of  system $\left(\mathrm{S}\right).$ The proof is
mostly
based on finding pairs of sub-solutions and super solutions of 
system $\left(\mathrm{S}\right).$ To this end, first recall from Candito et al.  \cite[Lemma 2]{CaMaMo} that
 $z_i\in \mathcal{C}^{1,\gamma }(\overline{\Omega })$ for $i=1,2,$ and for some $\gamma\in(0,1)$ are the unique solutions of  the homogeneous Dirichlet
problem 
\begin{equation}\label{eq26}
\begin{cases}
-\Delta _{p_i}u=1 \quad\text{ in }\Omega,\\
\text{ \ }u=0\qquad\text{ on }\partial \Omega,
\end{cases}
\end{equation}%
which satisfies 
\begin{equation}\label{eq27}
|| z_i||_{\mathcal{C}^{1,\gamma }(\overline{\Omega })}\leq L\qquad\textrm{and}\qquad||\nabla z_i||_\infty\leq \hat{L}\text{,}  
\end{equation}%
\begin{equation}\label{eq28}
ld(x)\leq z_i\leq Ld(x)\text{\ in\ }\Omega \text{,\ \ \ \ }\frac{%
\partial z_i}{\partial \eta }<0\text{\ on\ }\partial \Omega \text{,}
\end{equation}%
for certain constants $\hat{L}, l,$ and $ L.$ Moreover,  
by the
 Minty-Browder Theorem (see Brezis \cite{B})
combined with the Lieberman regularity Theorem \cite{L}, we know that the Dirichlet
problem 
\begin{equation}  \label{eq29}
-\Delta_{{p}_i}u=
\begin{cases}
1 & \text{if }x\in \Omega \backslash \overline{\Omega }_{\tau }, \\ 
-1 & \text{otherwise},%
\end{cases} \quad u=0\text{ on }\partial \Omega , 
\end{equation}
has a unique solution denoted by  $z_{i,\tau }\in 
\mathcal{C}^{1,\gamma }(\overline{\Omega })$ for a given 
 $0<\tau <\mathrm{diam}(\Omega )$, satisfying
 \begin{equation}\label{eeq30}
 z_{i,\tau }\leq z_i\text{ \ in }\Omega  
\end{equation}%
\begin{equation}\label{eq30}
\frac{\partial z_{i,\tau }}{\partial \eta }<\frac{1}{2}\frac{\partial z_i}{%
\partial \eta }<0\text{ on }\partial \Omega \text{ \ and \ }z_{i,\tau }\geq 
\frac{1}{2}\text{\thinspace }z_i\text{ \ in }\Omega.
\end{equation}%
Now, for a given $\tau>0$, we define  
\begin{equation}\label{eq31}
\underline{u}_1:=\tau ^{\frac{1}{p_1}}z_{1,\tau }^{\omega_1 }-\tau,\qquad\qquad
\underline{u}_2:=\tau ^{\frac{1}{p_2}}z_{2,\tau }^{\omega_2 }-\tau,
\end{equation}
\begin{equation}\label{eq32}
\overline{u}_1:=\tau^{-p_1}z_1^{\overline{\omega}_1}-\tau, \qquad\qquad
\overline{u}_2:=\tau^{-p_2}z_2^{\overline{\omega}_2}-\tau,
\end{equation}%
where 
\begin{equation} \label{eq33}
\frac{\omega_i -1}{\omega_i }>\frac{1}{p_i-1}>\frac{\overline{\omega}_i-1}{%
\overline{\omega}_i}\text{ \ with }\omega_i >\overline{\omega}_i>1  
\end{equation}%
and%
\begin{equation}\label{eq34}
\overline{\omega_i }<1+p_i\Big(1-\frac{\max \{\alpha_i ,\beta_i \}}{p_i-1}\Big).  
\end{equation}%
According to equations \eqref{eq27}-\eqref{eq28}, ones has 
\begin{equation}\label{eq35}
\overline{u}_1\leq \tau^{-p_1}(Ld)^{\overline{\omega}_1},\qquad\qquad\overline{u}_2\leq \tau^{-p_2}(Ld)^{\overline{\omega}_2}
\end{equation}
\begin{equation}\label{eq36}
||\nabla \overline{u}_1||_{\infty }\leq \tau^{-p_1}\hat{L}_1, \qquad\qquad ||\nabla \overline{u}_2||_{\infty }\leq \tau^{-p_2}\hat{L}_2,
\end{equation}%
with $\hat{L}_i:=\overline{\omega}_iL^{\overline{\omega}_i}$ for $i=1,2.$ Furthermore, we have 
\begin{equation}\label{eq37}
\begin{cases}
\displaystyle\frac{\partial \overline{u}_1}{\partial \eta }=\tau ^{-p_1}\frac{\partial (z_1^{\overline{\omega}_1})}{\partial \eta }=\tau^{-p_1}\overline{\omega}_1z_1^{\overline{\omega}_1-1}
\frac{\partial z_1}{\partial \eta }=0,\\
\\
\displaystyle\frac{\partial \overline{u}_2}{\partial \eta }=\tau ^{-p_2}\frac{\partial (z_2^{\overline{\omega}_2})}{\partial \eta }=\tau^{-p_2}\overline{\omega}_2z_2^{\overline{\omega}_2-1}
\frac{\partial z_2}{\partial \eta }=0,
\end{cases} 
\end{equation}
 on $\partial \Omega,$ since $z_i$  is a solution of the Dirichlet problem  \eqref{eq26} and also for 
  $\omega_i ,\overline{\omega}_i>1,$  $i=1,2.$
  
Now, we shall prove the following result.
\begin{lemma}\label{lem8}
For a sufficiently small $\tau >0,$ we have $\underline{u}_1\leq \overline{u}_1$ and $\underline{u}_2\leq \overline{u}_2.$
\end{lemma}
\begin{proof}
At first, we show that $\underline{u}_1\leq \overline{u}_1$ in 
$\Omega.$ Taking a direct computation, we obtain
\begin{equation*}
\begin{array}{l}
\overline{u}_1(x)-\underline{u}_1(x)=\left( \tau ^{-p_1}z_1^{\bar{\omega_1}}-\tau
\right) -\left( \tau ^{\frac{1}{p_1}}z_{i,\tau }^{\omega_i }-\tau \right) \\ 
\geq \tau ^{-p_1}z_1^{\bar{\omega_1}}-\tau ^{\frac{1}{p_1}}z_1^{\omega_1 }=z^{\omega_1
}(\tau ^{-p_1}z_1^{\bar{\omega}_1-\omega_1 }-\tau ^{\frac{1}{p_1}}) \\ 
\geq z_1^{\omega_1 }(\tau ^{-p_1}(cd(x))^{\bar{\omega}_1-\omega_1 }-\tau ^{\frac{1}{p_1}})\geq 0,
\end{array}%
\end{equation*}%
since $\omega_1 >\overline{\omega }_1$ and $z_{1,\tau }\leq z_1$ for every small enough
$\tau <\mathrm{diam}(\Omega ).$ 
Therefore, $\underline{u}_1\leq \overline{u}_1$ in 
$\Omega.$ Finally, using a similar argument as above  we can obtain that $\underline{u}_2\leq \overline{u}_2$
 in 
$\Omega.$ 
\end{proof}

\section{Proofs of the Main Results}\label{s6}

\begin{proof}[Proof of Theorem~\ref{T1}]
First, we claim that equation \eqref{defsupersol} is  satisfied by the pair of functions
$(\overline{u}_1,\overline{u}_2)$ given by equation \eqref{eq32}. To see this, pick $(u_1,u_2)\in W^{1,p_1}(\Omega )\times W^{1,p_2}(\Omega )$ within $[\underline{u}_1,\overline{u}_1]\times [\underline{u}_2,\overline{u}_2]$ such that $\underline{
u}_2\leq u_2\leq \overline{u}_2$, $\underline{
u}_1\leq u_1\leq \overline{u}_1$. Now, in view of  condition $({\bf H_1}),$ combined with equations \eqref{eq35} and \eqref{eq36}, ones has
\begin{equation}\label{eq38}
\begin{split}
\Big|f_1(. ,\overline{u}_1,u_2,\nabla \overline{u}_1,\nabla u_2)\Big|& \leq 
M_1(1+|\overline{u}_1|^{\alpha_1 })(1+|\nabla\overline{u}_1 |^{\beta_1 })\\
&\leq M_1(1+(\tau ^{-p_1}(Ld)^{\overline{\omega }_1})^{\alpha_1 })(1+(\tau ^{-p_1}
\hat{L})^{\beta_1 })\\
&\leq 2M_1(C_1C_2)^{\alpha_1 \beta_1 }\tau ^{-p_1\max \{\alpha_1 ,\beta_1 \}}\\
&\leq C\tau ^{-p_1\max \{\alpha_1 ,\beta_1 \}},
\end{split}
\end{equation}%
where $C:=2M_1(C_1C_2)^{\alpha_1 \beta_1 }$ and $\tau >0$ is  small enough. Using the same argument as in equation \eqref{eq38}, we obtain

\begin{equation}\label{eq39}
\Big|f_2(\cdot ,u_1,\overline{u}_2,\nabla u_1,\nabla \overline{u}_2)\Big|\leq  C'\tau ^{-p_2\max \{\alpha_2 ,\beta_2 \}}\text{,}
\end{equation}%
for some constant $C'>0$ and for $\tau >0$  small enough.
Now, in view of equations \eqref{eq32}
and
\eqref{eq33}, we have
\begin{equation}\label{eq40}
-\Delta _{p_1}\overline{u}_1+\frac{|\nabla \overline{u}_1|^{p_1}}{\overline{u}_1%
+\delta_1 }=\tau ^{-p_1(p_1-1)}\Big(-\Delta _{p_1}z_1^{\bar{\omega}_1}+\frac{|\nabla z_1^{%
\bar{\omega}_1}|^{p}_1}{z_1^{\bar{\omega}_1}}\Big).
\end{equation}
On the other hand, with a direct
computation,  we can get
\begin{eqnarray}\label{eq41}
-\Delta _{p_1}z_1^{\bar{\omega}_1}+ \frac{|\nabla z_1^{\bar{\omega}_1}|^{p_1}}{z_1^{%
\bar{\omega}_1}}&=&\bar{\omega}_1^{p_1-1}\Big( 1-(\bar{\omega}_1-1)( p_1-1) 
\frac{|\nabla z_1|^{p_1}}{z_1}\Big) z_1^{(\bar{\omega_1}-1)(p_1-1)}\nonumber\\
&+&\bar{\omega_1}^{p_1}%
\frac{z_1^{(\bar{\omega}_1-1)p_1}|\nabla z_1|^{p_1}}{z_1^{\bar{\omega}_1}}\nonumber \\ 
&=&\bar{\omega}^{p_1-1}\Big( 1-(\bar{\omega}_1-1)( p_1-1) \frac{|\nabla z_1|^{p_1}}{z_1}\Big) z_1^{(\bar{\omega}_1-1)(p_1-1)}\nonumber\\
&+&\bar{\omega}_1^{p_1}z_1^{(\bar{\omega}_1-1)(p_1-1)}\frac{z_1^{\bar{\omega}_1-1}|\nabla z_1|^{p_1}}{z_1^{\bar{\omega}_1}} \nonumber\\ 
&=&\bar{\omega}_1^{p_1-1}\Big[ 1+\bar{\omega}_1(1-\frac{(\bar{\omega}_1-1)(p_1-1) }{\bar{\omega}_1})\frac{|\nabla z_1|^{p_1}}{z_1}\Big] z_1^{(\bar{\omega}_1-1)(p_1-1)}.\nonumber\\
& &
\end{eqnarray}
Invoking equations \eqref{eq40} and  \eqref{eq41}, it follows that
\begin{eqnarray}\label{eq42}
-\Delta _{p_1}\overline{u}_1+\frac{|\nabla \overline{u}_1|^{p_1}}{\overline{u}_1%
+\delta_1 } 
=\tau ^{-p_1(p_1-1)}\bar{\omega}_1^{p_1-1}\Big[ 1+\bar{\omega}_1(1-\frac{(\bar{%
\omega}_1-1)( p_1-1) }{\bar{\omega}_1})\frac{|\nabla z_1|^{p_1}}{z_1}\Big]
z_1^{(\bar{\omega}_1-1)(p_1-1)} \nonumber\\ \nonumber
\geq\displaystyle \tau ^{-p_1(p_1-1)}\bar{\omega}_1^{p_1-1}\left\{ 
\begin{array}{ll}\displaystyle
z_1^{(\bar{\omega}_1-1)(p_1-1)} & \text{in}\;\Omega \backslash \overline{\Omega }%
_{\tau }, \\ 
\\
\bar{\omega}_1(1-\frac{(\bar{\omega}_1-1)( p_1-1) }{\bar{\omega}_1})z_1^{(%
\bar{\omega}_1-1)(p_1-1)-1}|\nabla z_1|^{p_1} & \text{in}\;\Omega _{\tau}\text{.}
\end{array}%
\right.
\end{eqnarray}
Moreover, ones has after using equation \eqref{eq34} and decreasing $\tau $ if necessary,
\begin{eqnarray}\label{eq43}
\displaystyle\tau ^{-p_1(p_1-1)}\bar{\omega}_1^{p_1-1}z_1^{(\bar{\omega}_1-1)(p_1-1)} 
&\geq& \tau ^{-p_1(p_1-1)}\bar{\omega}^{p_1-1}(c^{-1}d(x))^{(\bar{\omega}_1-1)(p_1-1)}\nonumber
\\ 
&\geq& \tau ^{-p_1(p_1-1)}\bar{\omega}_1^{p_1-1}(c^{-1}\tau )^{(\bar{\omega}%
-1)(p_1-1)}\nonumber\\ 
&=&\tau ^{(\bar{\omega}_1-1-p_1)(p_1-1)}\bar{\omega}_1^{p_1-1}c^{-(\bar{\omega}_1
-1)(p_1-1)}\nonumber \\ 
&\geq& \tau ^{-p_1\max \{\alpha_1 ,\beta_1 \}}\text{ \ in }\Omega \backslash 
\overline{\Omega }_{\tau }.
\end{eqnarray}
Finally, combining equations \eqref{eq38} and \eqref{eq43}, we obtain
\begin{equation}\label{eq44}
-\Delta _{p_1}\overline{u}_1+\frac{|\nabla \overline{u}_1|^{p_1}}{\overline{u}_1
+\delta_1 }\geq f_1(\cdot ,\overline{u}_1,u_2,\nabla \overline{u}_1,\nabla u_2)\text{\ \ in\ }%
\Omega \backslash \overline{\Omega }_{\tau }\text{.}
\end{equation}
Now, pick any  $x\in {\Omega }_{\tau }.$ By virtue of equations \eqref{eq28} and  \eqref{eq33}, we can find a constant $\beta>0$ such that 
\begin{equation*}
(1-\frac{(\bar{\omega}_1-1)\left( p_1-1\right) }{\bar{\omega}_1})|\nabla z_1|>\beta\text{ \ in }{\Omega }_{\tau }.
\end{equation*}%
By equations \eqref{eq28} and \eqref{eq33}, one has
\begin{eqnarray*}
&&\displaystyle\tau^{-p_1(p_1-1)}\bar{\omega}_1^{p_1}\Big(1-\frac{(\bar{\omega}_1-1)\left( p_1-1\right) 
}{\bar{\omega}_1}\Big)z_1^{(\bar{\omega}_1-1)(p_1-1)-1}|\nabla z_1{|}^{p_1}\nonumber\\ 
&\geq&\tau^{-p_1(p_1-1)}\bar{\omega}_1^{p_1}(Ld(x))^{(\bar{\omega}_1-1)(p_1-1)-1}\bar{%
\mu}^{p} \\ 
&\geq&\tau^{-p_1(p_1-1)}\bar{\omega}_1^{p_1}(L\tau )^{(\bar{\omega}_1-1)(p_1-1)-1}\beta^{p_1} \\ 
&\geq& \tau^{-p_1\max \{\alpha_1 ,\beta_1 \}}\text{ \ in }\Omega _{\tau }.
\end{eqnarray*}%
Therefore, in view of the above equation and for $\tau >0$ sufficiently  small, we obtain 
\begin{equation}\label{eq45}
-\Delta _{p_1}\overline{u}_1+\frac{|\nabla \overline{u}_1|^{p_1}}{\overline{u}_1
+\delta_1 }\geq f_1(\cdot ,\overline{u}_1,u_2,\nabla \overline{u}_1,\nabla u_2)\text{\ \ in\ }%
\Omega _{\tau}\text{.}
\end{equation}%
Combining equations \eqref{eq44}
and
 \eqref{eq45},  we get
\begin{equation}\label{eq46}
-\Delta _{p_1}\overline{u}_1+\frac{|\nabla \overline{u}_1|^{p_1}}{\overline{u}_1
+\delta_1 }\geq f_1(\cdot ,\overline{u}_1,u_2,\nabla \overline{u}_1,\nabla u_2)\text{\ \ in\ }%
\Omega\text{.}
\end{equation}
A similar argument yields
\begin{equation}\label{eq47}
-\Delta _{p_2}\overline{u}_2+\frac{|\nabla \overline{u}_2|^{p_2}}{\overline{u}_2
+\delta_2 }\geq f_2(\cdot ,u_1,\overline{u}_2,\nabla u_1,\nabla\overline{u}_2)\text{\ \ in\ }%
\Omega\text{.}
\end{equation}
Now, test equation \eqref{eq46} and equation \eqref{eq47} with $(\varphi_1,\varphi_2) \in W_{b}^{1,p_1}(\Omega )\times W_{b}^{1,p_2}(\Omega ),$  $\varphi \geq 0$ a.e.
in $\Omega $, and equation \eqref{eq37} yield
\begin{eqnarray*}
\int_{\Omega }|\nabla \overline{u}_1|^{p_1-2}\nabla \overline{u}_1\nabla \varphi_1 \,%
\mathrm{d}x&+&\int_{\Omega }\frac{|\nabla \overline{u}_1|^{p_1}}{\overline{u}_1
+\delta_1}\varphi_1 \mathrm{d}x-\left\langle \frac{\partial \overline{u}}{%
\partial \eta _{p_1}},\gamma _{0}(\varphi_1 )\right\rangle _{\partial \Omega }
\\ 
&\geq&\displaystyle \int_{\Omega }f_1(\cdot ,\overline{u}_1,u_2,\nabla \overline{u}_1,\nabla u_2)\varphi_1 \text{%
\thinspace \textrm{d}}x\text{,}%
\end{eqnarray*}
\begin{eqnarray*}
\int_{\Omega }|\nabla \overline{u}_2|^{p_2-2}\nabla \overline{u}_2\nabla \varphi_2 \,%
\mathrm{d}x&+&\int_{\Omega }\frac{|\nabla \overline{u}_2|^{p_2}}{\overline{u}_2
+\delta_2}\varphi_2 \mathrm{d}x-\left\langle \frac{\partial \overline{u}}{%
\partial \eta _{p_2}},\gamma _{0}(\varphi_2)\right\rangle _{\partial \Omega }
\\ 
&\geq&\displaystyle \int_{\Omega }f_2(\cdot ,u_1,\overline{u}_2,\nabla u_1,\nabla\overline{u}_2)\varphi_2 \text{%
\thinspace \textrm{d}}x\text{,}%
\end{eqnarray*}
where $\gamma _{0}$ is the trace operator on $\partial \Omega $, 
\begin{equation}
\frac{\partial w}{\partial \eta _{p_i}}:=|\nabla w|^{p_i-2}\frac{\partial w}{%
\partial \eta },\quad \hbox{for every} \,w\in W^{1,p_i}(\Omega )\cap C^{1}(\overline{%
\Omega }),  \label{conode}
\end{equation}%
and  $\left\langle \cdot ,\cdot \right\rangle _{\partial \Omega }$ is 
the duality brackets for the pair 
\begin{equation*}
(W^{1/p_i^{\prime },p_i}(\partial \Omega ),W^{-1/p_i^{\prime },p_i^{\prime
}}(\partial \Omega )). 
\end{equation*}
The proof of the claim is now completed.

Next, we show  that equation \eqref{defsubsol} is  satisfied by the pair of functions
$(\underline{u}_1,\underline{u}_2)$ given by equation \eqref{eq31}. A direct
computation yields
\begin{eqnarray*}
-\Delta _{p_1}z_{1,\tau }^{\omega_1 }+\frac{|\nabla z{1,\tau }^{\omega_1 }|^{p_1}}{%
z_{1,\tau }^{\omega_1 }}&=&\omega_1^{p_1-1}\left( 1-(\omega_1 -1)\left( p_1-1\right) 
\frac{|\nabla z_{1,\tau }|^{p_1}}{z_{1,\tau }}\right) z_{1,\tau }^{(\omega_1
-1)(p_1-1)}\\&+&\omega_1^{p_1}\frac{z_{1,\tau }^{(\omega_1 -1)p_1}|\nabla z_{1,\tau }|^{p_1}%
}{z_{1,\tau }^{\omega_1}} \\ 
&=&\omega_1^{p_1-1}\Big[ 1+\omega_1 (1-\frac{(\omega_1 -1)\left( p_1-1\right) }{\omega_1 
})\frac{|\nabla z_{1,\tau }|^{p_1}}{z_{1,\tau }}\Big] z_{1,\tau }^{(\omega_1
-1)(p_1-1)},
\end{eqnarray*}
in $ \Omega \backslash \overline{\Omega }_{\tau }.$ Similarly,  in $\Omega _{\tau }$ it follows that
\begin{eqnarray*}
-\Delta _{p_1}z_{1,\tau}^{\omega_1 }+\frac{|\nabla z_{1,\tau }^{\omega_1 }|^{p_1}}{%
z_{1,\tau }^{\omega_1 }}=\omega_1^{p_1-1}\Big[ -1+\omega_1 (1-\frac{(\omega_1
-1)\left( p_1-1\right) }{\omega })\frac{|\nabla z_{1,\tau }|^{p_1}}{z_{1,\tau}}%
\Big] z_{1,\tau }^{(\omega_1 -1)(p_1-1)}.
\end{eqnarray*}
In fact,  by  equations \eqref{eq31} and   \eqref{eq33}, ones has 
\begin{equation}\label{eq48}
-\Delta _{p_1}\underline{u}_1+\frac{|\nabla \underline{u}_1|^{p_1}}{\underline{u}_1
+\delta_1 }=\tau ^{\frac{1}{p_1^{\prime }}}(-\Delta _{p_1}z_{1,\tau }^{\omega_1}+%
\frac{|\nabla z_{1,\tau}^{\omega_1 }|^{p_1}}{z_{1,\tau }^{\omega_1 }})\leq
\left\{ 
\begin{array}{ll}
\tau ^{\frac{1}{p_1^{\prime }}}\omega_1^{p_1-1}z_{1,\tau}^{(\omega_1 -1)(p_1-1)} & 
\text{in }\Omega \backslash \overline{\Omega }_{\tau} \\ 
\\
0 & \text{in }\Omega_{\tau }.%
\end{array}%
\right.
\end{equation}%
In view  of equations \eqref{eq27}-\eqref{eeq30}, after having chosen an appropriate constant 
$m_1 $ in $({\bf H_2}),$ 
\begin{equation}\label{eq49}
\displaystyle m_1>\tau^{\frac{1}{p_1^{\prime }}}\omega_1^{p_1-1}L^{(\omega_1 -1)(p_1-1)}\text{ \
for }\tau >0\text{ small enough.}
\end{equation}
Combining equations \eqref{eq48}
and \eqref{eq49}, we arrive at
\begin{equation}\label{eq50}
\displaystyle-\Delta _{p_1}\underline{u}_1+\frac{|\nabla \underline{u}_1|^{p_1}}{\underline{u}_1
+\delta_1 }\leq f_1(\cdot ,\underline{u}_1,u_2,\nabla \underline{u}_1,\nabla u_2).
\end{equation}
Using a similar argument as above, we obtain
\begin{equation}\label{eq51}
\displaystyle-\Delta _{p_2}\underline{u}_2+\frac{|\nabla \underline{u}_2|^{p_2}}{\underline{u}_2
+\delta_2 }\leq f_2(\cdot ,u_1,\underline{u}_2,\nabla u_1, \underline{u}_2).    
\end{equation}
Finally, test equations \eqref{eq50} and \eqref{eq51} with $(\varphi_1, \varphi_2) \in W_{b}^{1,p_1}(\Omega )\times W_{b}^{1,p_2}(\Omega ),$ where $\varphi_1, \varphi_2
\geq 0$ a.e. in $\Omega ,$ and use equation \eqref{eq37} and the Green formula 
\cite{CF}, to obtain
\begin{eqnarray*}
\int_{\Omega }|\nabla \underline{u}_1|^{p_1-2}\nabla \underline{u}_1\nabla \varphi_1
\,\mathrm{d}x&+&\int_{\Omega }\frac{|\nabla \underline{u}_1|^{p_1}}{\underline{u}_1
+\delta_1 }\varphi_1 \,\mathrm{d}x \\ 
&\leq& \int_{\Omega }|\nabla \underline{u}_1|^{p_1-2}\nabla \underline{u}_1\nabla
\varphi_1 \,\mathrm{d}x-\left\langle \frac{\partial \underline{u}_1}{\partial
\eta _{p_1}},\gamma _{0}(\varphi_1 )\right\rangle _{\partial \Omega
}+\int_{\Omega }\frac{|\nabla \underline{u}_1|^{p_1}}{\underline{u}_1+\delta_1 }
\varphi_1 \,\mathrm{d}x \\ 
&=&\int_{\Omega }(-\Delta _{p_1}\underline{u}_1\ +\frac{|\nabla \underline{u}_1|^{p_1}%
}{\underline{u}_1+\delta_1 })\varphi_1 \,\mathrm{d}x\\
&\leq& \int_{\Omega }f_1(\cdot ,%
\underline{u}_1,u_2,\nabla \underline{u}_1,\nabla u_2)\varphi_1 \,\mathrm{d}x,%
\end{eqnarray*}
\begin{eqnarray*}
\int_{\Omega }|\nabla \underline{u}_2|^{p_2-2}\nabla \underline{u}_2\nabla \varphi_2
\,\mathrm{d}x&+&\int_{\Omega }\frac{|\nabla \underline{u}_2|^{p_2}}{\underline{u}_2
+\delta_2 }\varphi_2 \,\mathrm{d}x \\ 
&\leq& \int_{\Omega }|\nabla \underline{u}_2|^{p_2-2}\nabla \underline{u}_2\nabla
\varphi_2 \,\mathrm{d}x-\left\langle \frac{\partial \underline{u}_2}{\partial
\eta_{p_2}},\gamma _{0}(\varphi_2 )\right\rangle _{\partial \Omega
}+\int_{\Omega }\frac{|\nabla \underline{u}_2|^{p_2}}{\underline{u}_2+\delta_2 }
\varphi_2 \,\mathrm{d}x \\ 
&=&\int_{\Omega }(-\Delta _{p_2}\underline{u}_2 +\frac{|\nabla \underline{u}_2|^{p_2}%
}{\underline{u}_2+\delta_2})\varphi_2 \,\mathrm{d}x\\
&\leq& \int_{\Omega }f_2(\cdot, u_1
\underline{u}_2,\nabla u_1, \nabla \underline{u}_2)\varphi_2 \,\mathrm{d}x,%
\end{eqnarray*}
since $\gamma _{0}(\varphi_1 ),\,\gamma _{0}(\varphi_2 )\geq 0,$ whenever $(\varphi_1,\varphi_2) \in W^{1,p_1}(\Omega )\times W^{1,p_2}(\Omega )$ (for more details see Carl et al.  \cite[p. 35]{CLM}).

Consequently, $(\underline{u}_1,\underline{u}_2)$ and $(\overline{u}_1,\overline{u}_2)$ satisfy equations \eqref{c3}
and
\eqref{c4}. Therefore, we can apply  Theorem \ref{thmsubsuper}, and we obtain the existence of  a solution $(u_{0},u'_{0})\in \mathcal{C}^{1,\gamma }(\overline{\Omega })\times\mathcal{C}^{1,\gamma }(\overline{\Omega })$  of system $\left( \mathrm{S}\right),$ satisfying
\begin{equation}\label{eq52}
\underline{u}_1\leq u_{0}\leq \overline{u}_1,\qquad\qquad \underline{u}_2\leq u'_{0}\leq \overline{u}_2.  
\end{equation}%
Furthermore, $(u_{0},u'_{0})$ is nodal solution. Indeed, combining equations \eqref{eq28},
 \eqref{eq31},  and \eqref{eq32}, we arrive at
\begin{equation*}
\overline{u}_1=\tau^{-p_1}z^{\overline{\omega }_1}-\tau \leq \tau^{-p_1}(Ld(x))^{\overline{\omega }_1}-\tau ,
\end{equation*}
\begin{equation*}
\overline{u}_2=\tau^{-p_2}z^{\overline{\omega }_2}-\tau \leq \tau^{-p_2}(Ld(x))^{\overline{\omega }_2}-\tau ,
\end{equation*}
which implies that
\begin{equation}\label{eq53}
\max\{\overline{u}_1(x),\overline{u}_2(x)\}<0,\;\;\text{provided that}\;\;d(x)<L^{-1}\tau^{\frac{p_i+1}{%
\overline{\omega }_i}},
\end{equation}%
for $i=1,2.$ Combining  equations \eqref{eq28},
 \eqref{eq31},
 and
  \eqref{eq32}, yields
\begin{equation*}
\underline{u}_1=\tau^{\frac{1}{p_1}}z_{1,\tau}^{\omega_1 }-\tau \geq \tau^{\frac{1}{p_1}}(ld(x))^{\omega_1}-\tau ,
\end{equation*}
\begin{equation*}
\underline{u}_2=\tau^{\frac{1}{p_2}}z_{2,\tau}^{\omega_2 }-\tau \geq \tau^{\frac{1}{p_2}}(ld(x))^{\omega_2}-\tau,
\end{equation*}
hence
\begin{equation}\label{eq54}
\min\{\underline{u}_1(x),\underline{u}_2(x)\}>0,\;\;\text{when}\;\;d(x)>l\tau^{\frac{1}{\omega_i
p_i^{\prime }}}\text{.}  
\end{equation}
 for $i=1,2.$ 
 The conclusion now follows from equations \eqref{eq53}
 and \eqref{eq54}. 
 This completes the proof of Theorem~\ref{T1}.
\end{proof}

\begin{proof}[Proof of Theorem~\ref{Thm2}]
    First, using the same notation as in equations \eqref{eq31} and \eqref{eq32} and applying the same argument as
    in the proof of 
    Theorem \ref{T1}, we can
     ensure that $(\underline{u}_{+},\underline{u}^{+})$
    and $(\overline{u}_{+},\overline{u}^{+})$ satisfy the equations \eqref{c3}  and \eqref{c4}. Therefore, invoking Theorem \ref{thmsubsuper}, we obtain the existence of a solution $(u_{+},u^{+})\in \mathcal{C}^{1,\gamma }(\overline{\Omega })\times \mathcal{C}^{1,\gamma }(\overline{\Omega })$ with the following properties:
    \begin{equation*}
        u_{0}\leq u_{+}\leq \overline{u}_{+}\quad \textrm{and}\;\;
u_{+}\geq 0 \;\;\textrm{on}\; \Omega,
    \end{equation*}
    \begin{equation*}
        u'_{0}\leq u^{+}\leq \overline{u}^{+}\quad \textrm{and}\;\;
u^{+}\geq 0 \;\;\textrm{on}\; \Omega.
    \end{equation*}
 Finally,  using equation \eqref{eq32}, we can easily deduce that $u_{+}(x)$ and $u^{+}(x)$ are zero, when $%
d(x)\rightarrow 0$.
This completes the proof of Theorem~\ref{Thm2}.
\end{proof}

\subsection*{Acknowledgements}
 The authors gratefully acknowledge comments and suggestions by the reviewers.
 Repov\v{s} was supported by the Slovenian Research and Innovation Agency
 grants number P1-0292, J1-4031, J1-4001, N1-0278, N1-0114, and N1-0083.

\end{document}